\documentclass[12pt]{article}
\usepackage[margin=1in,headheight=15pt]{geometry}
\usepackage[T1]{fontenc}
\usepackage{iftex}
\IfFileExists{gfsartemisia-euler.sty}
  {\usepackage{gfsartemisia-euler}
   \ifPDFTeX\else
     \usepackage[no-math]{fontspec}
     \setmainfont{GFSArtemisia.otf}[BoldFont=GFSArtemisiaBold.otf,
       ItalicFont=GFSArtemisiaIt.otf,BoldItalicFont=GFSArtemisiaBoldIt.otf]
   \fi
   \AtBeginDocument{\SetMathAlphabet{\mathbf}{normal}{OT1}{artemisiaeuler}{b}{n}%
                    \SetMathAlphabet{\mathbf}{bold}{OT1}{artemisiaeuler}{b}{n}}}
  {\usepackage{lmodern}}
\usepackage{amsmath,amssymb,amsthm,mathtools,enumitem,braket,authblk}
\usepackage[cal=euler,frak=euler]{mathalfa}
\usepackage{booktabs,aliascnt,tikz-cd,fancyhdr,xcolor}
\usepackage{microtype}
\usepackage{quiver}
\usetikzlibrary{nfold}
\usepackage[backend=biber,style=alphabetic,maxbibnames=99,giveninits=true]{biblatex}
\usepackage{hyperref}
\usepackage[nameinlink,noabbrev]{cleveref}

\definecolor{Unipv}{RGB}{165, 56, 77}
\definecolor{bluegreen}{RGB}{8, 143, 143}

\hypersetup{
    colorlinks=true,
    linkcolor=Unipv,
    citecolor=bluegreen,
    urlcolor=blue
}

\theoremstyle{plain}
\newtheorem{theorem}{Theorem}[section]

\newaliascnt{proposition}{theorem}
\newtheorem{proposition}[proposition]{Proposition}
\aliascntresetthe{proposition}

\newaliascnt{lemma}{theorem}
\newtheorem{lemma}[lemma]{Lemma}
\aliascntresetthe{lemma}

\newaliascnt{corollary}{theorem}
\newtheorem{corollary}[corollary]{Corollary}
\aliascntresetthe{corollary}

\theoremstyle{definition}

\newaliascnt{definition}{theorem}
\newtheorem{definition}[definition]{Definition}
\aliascntresetthe{definition}

\newaliascnt{example}{theorem}
\newtheorem{example}[example]{Example}
\aliascntresetthe{example}

\theoremstyle{remark}

\newaliascnt{remark}{theorem}
\newtheorem{remark}[remark]{Remark}
\aliascntresetthe{remark}

\crefname{theorem}{theorem}{theorems}
\Crefname{theorem}{Theorem}{Theorems}

\crefname{proposition}{proposition}{propositions}
\Crefname{proposition}{Proposition}{Propositions}

\crefname{lemma}{lemma}{lemmas}
\Crefname{lemma}{Lemma}{Lemmas}

\crefname{corollary}{corollary}{corollaries}
\Crefname{corollary}{Corollary}{Corollaries}

\crefname{definition}{definition}{definitions}
\Crefname{definition}{Definition}{Definitions}

\crefname{example}{example}{examples}
\Crefname{example}{Example}{Examples}

\crefname{remark}{remark}{remarks}
\Crefname{remark}{Remark}{Remarks}

\newcommand{\Man}{\mathbf{Man}}
\newcommand{\St}{\mathbf{St}}

\newcommand{\Hom}{\operatorname{Hom}}

\newcommand{\Mor}{\operatorname{Mor}}
\newcommand{\Aut}{\operatorname{Aut}}

\newcommand{\Stab}{\operatorname{St}}

\newcommand{\id}{\operatorname{id}}
\newcommand{\pt}{\mathrm{pt}}

\newcommand{\cX}{\mathcal X}
\newcommand{\cY}{\mathcal Y}
\newcommand{\cM}{\mathcal M}
\newcommand{\cN}{\mathcal N}
\newcommand{\cW}{\mathcal W}
\newcommand{\cU}{\mathcal U}
\newcommand{\rX}{\mathsf X}
\newcommand{\rY}{\mathsf Y}
\newcommand{\BG}{\mathcal{B}G}

\newcommand{\st}{\mathrm{st}}
\newcommand{\rs}{\mathsf{s}}
\newcommand{\rt}{\mathsf{t}}
\newcommand{\rr}{\mathsf{r}}

\newcommand{\acts}{\cdot}

\newcommand{\Sone}{\mathbb{S}^1}
\newcommand{\C}{\mathbb C}
\newcommand{\HomotopyFixedLocus}{\mathrm{H}}
\newcommand{\FixedLocus}{\mathrm{F}}
\newcommand{\Forgetful}{\mathrm{U}}
\newcommand{\Trivial}{\mathrm{T}}
\newcommand{\qman}{M}

\title{On the Fixed Loci of a Differentiable Stack}
\author[1,2]{Leone Parasporo\thanks{\texttt{l.parasporo@campus.unimib.it}}}
\affil[1]{Università degli Studi di Pavia, Dipartimento di Matematica, Via Ferrata 1, 27100 Pavia, Italy}
\affil[2]{Università degli Studi di Milano-Bicocca, Dipartimento di Matematica e Applicazioni, Via R. Cozzi 55, 20126 Milan, Italy}

\date{November 26, 2026}

\begin{document}

\maketitle

\begin{abstract}
Let $G$ be a compact Lie group acting smoothly on a separated differentiable stack $\cM$.
We prove that the homotopy fixed locus $\cM^{hG}$ is a separated differentiable stack with quasi-representable evaluation map to $\cM$.
Its local sectors correspond to conjugacy classes of splittings of the isotropy extension.
If $G$ is connected and isotropy dimension is locally constant on the pointwise fixed locus, the stabilizer-theoretic fixed locus $\cM^G$ is regular and differentiable.
The comparison map $\cM^{hG}\to\cM^G$ is then a representable submersion. It is proper precisely when each fixed point has finitely many splitting classes.
\end{abstract}

\clearpage

\section*{Introduction}

For a smooth action of a Lie group $G$ on a manifold $M$, the phrase ``fixed locus'' has an unambiguous meaning:
\[
M^G=\Set{m\in M|g\cdot m=m \ \forall\,g\in G}.
\]
However, for a smooth $G$-action on a differentiable stack $\cM$, there is an immediate ambiguity.  A point $x$ of a stack need not be fixed strictly by $g$; it may be fixed only up to an isomorphism
\[
g\cdot x\simeq x.
\]
Moreover, one can either ask merely for the existence of such isomorphisms, locally after every change of test manifold and for each smooth section $g$, or one can include the isomorphisms as part of the datum and impose the group law coherently.  These two possibilities lead to two different objects.

The first is the \emph{stabilizer-theoretic fixed locus}. In the algebraic and derived setting, Aranha, Khan, Latyntsev, Park and Ravi require the canonical stabilizer homomorphism to be locally surjective. Equivalently, the $G$-stabilizer of $x$ is all of $G$ \cite[Definition~A.9 and Remark~A.10]{aranha2025virtual}. The corresponding condition on smooth test families defines a full substack
\[
\cM^{G}\subseteq \cM.
\]

The second is the \emph{homotopy fixed locus} $\cM^{hG}$. Its algebraic counterpart is Romagny's fixed point stack, which represents equivariant morphisms from a stack with trivial action \cite[Definition~2.3 and Proposition~2.5]{romagny2005group}. Aranha et al.\ express the coherent choice of isomorphisms as a group-theoretic section of the extended stabilizer projection \cite[Definition~A.14 and Remarks~A.16--A.17]{aranha2025virtual}. For the homotopy-theoretic origins of the terminology, see for example \cite{carlsson1992survey,thomason1983homotopy}.

For compact $G$ and separated $\cM$, we describe $\cM^{hG}$ locally by applying the quotient calculation of \Cref{thm:homotopy-fixed-quotient} to a slice that preserves the classifying map $[\cM/G]\to\BG$ (\Cref{lem:relative-slice}). This proves differentiability and separatedness of $\cM^{hG}$ and gives its isotropy groups (\Cref{thm:hfp-differentiability}). Under the hypotheses of \Cref{thm:comparison}, we determine the image of $\cM^{hG}\to\cM^G$, compute its pullback to an arbitrary test manifold, and characterize properness in terms of splittings.

For a proper \'etale groupoid and a compact connected acting group, Jiang, Yang and Chen prove that the pointwise fixed set is embedded. Their proof uses a slice, the identity component of the extended isotropy group and a dimension count \cite[Definition~4.7, Lemma~4.8 and Theorem~4.9]{jiang2026equivariant}. We extend their argument to positive-dimensional isotropy under a local constancy assumption on the pointwise fixed locus, and verify the lifting condition for smooth families. In the orbifold case, \Cref{cor:fixed-locus-orbifold-stack} identifies the intrinsic stack presented by their fixed groupoid.

\Cref{sec:definitions} fixes our conventions for actions and presentations. \Cref{sec:fixed-locus} treats the stabilizer-theoretic fixed locus. \Cref{sec:homotopy-fixed-locus} computes the homotopy fixed locus for actions induced by group extensions, proves the general result for separated stacks, and compares the two loci. \Cref{sec:examples} illustrates the comparison for classifying stacks, circle actions and orbifolds. \Cref{appendix:translation} proves the correspondence between smooth stack actions and Kan fibrations over the classifying groupoid.

Throughout, $\Man$ is the site of Hausdorff, second-countable smooth manifolds of locally constant dimension without boundary, equipped with the \emph{open-cover topology}. We work in the $(2,1)$-category of stacks in groupoids on this site. A \emph{prestack} is a category fibered in groupoids whose isomorphism presheaves are sheaves. A \emph{stack} is a prestack with effective descent for objects.
An \emph{atlas} is a representable surjective submersion from a manifold. A \emph{differentiable stack} is a stack admitting an atlas \cite{behrend2011differentiable,heinloth2005notes,ginot2013introduction,del2012lie,bursztyn2020principal}.

A morphism $f:\cX\to\cY$ is \emph{representable} if $T\times_{\cY}\cX$ is represented by a manifold for every manifold $T\to\cY$. When the target is differentiable, we call $f$ \emph{quasi-representable} if this holds for an atlas of $\cY$, equivalently for every atlas \cite[Definition~2.6 and Proposition~2.7]{bursztyn2020principal}. For a representable morphism, submersiveness and properness are tested on every such $T$. A quasi-representable closed embedding means that every atlas pullback is a closed embedding of manifolds. All stack fiber products are $2$-fiber products. Their objects include a specified isomorphism between the two images, and their arrows intertwine this isomorphism. For stacks we use the terms \emph{isomorphism} and \emph{equivalence} interchangeably.

For Morita maps and Morita equivalences, we follow the conventions of \cite[\S\S~4.2--4.5]{del2012lie}.

We use $\cX,\cY$ for general (pre)stacks, $\cM,\cN$ for differentiable stacks and plain letters for manifolds. We write $\mathsf B G=(G\rightrightarrows\pt)$ for the one-object Lie groupoid, $\BG$ for its stack of principal right $G$-bundles, and $\mathcal B_G$ for the ordinary category of parameterized group elements in \Cref{rem:pseudofunctor}. In a Lie groupoid $\rX=(X_1\rightrightarrows X_0)$, $\rs,\rt$ are source and target, $1_x$ is the unit, and $ba=b\circ a$ when $\rs(b)=\rt(a)$. The inverse arrow is $a^{-1}$.

For a fixed Lie group $G$, we write
\[
\begin{aligned}
\Trivial&:\St\longrightarrow\St_G,\quad &
\Forgetful&:\St_G\longrightarrow\St,\\
\FixedLocus&:\St_G\longrightarrow\St,\quad &
\HomotopyFixedLocus&:\St_G\longrightarrow\St
\end{aligned}
\]
for the trivial-action, forgetful, fixed-locus and homotopy-fixed-locus functors, respectively. These are pseudofunctors of $(2,1)$-categories. We understand naturality and adjunctions bicategorically.

\subsection*{Acknowledgments}
The author thanks Professors Frank Neumann and Ping Xu for their helpful discussions, advice, and suggestions. The author is financially supported by the University of Pavia and is a member of GNSAGA–INdAM.

\clearpage

\tableofcontents

\section{Group actions and presentations}\label{sec:definitions}

\begin{definition}[\cite{romagny2005group, barbosa2021equivariant}]
	Let $G$ be a Lie group and $\mathcal{X}$ a prestack over $\Man$. A \emph{left $G$-action} on $\mathcal{X}$ consists of a morphism of prestacks $\mu:G\times\mathcal{X}\to\mathcal{X}$ and invertible $2$-cells $\iota$ and $\alpha$ such that
	\begin{enumerate}[label=(\roman*)]
		\item (identity) 
		\[\begin{tikzcd}
		{\mathcal{X}} & {G\times\mathcal{X}} \\
		& {\mathcal{X}}
		\arrow["{(e,\id_{\mathcal{X}})}", from=1-1, to=1-2]
		\arrow[""{name=0, anchor=center, inner sep=0}, "{\id_{\mathcal{X}}}"', from=1-1, to=2-2]
		\arrow["\mu", from=1-2, to=2-2]
		\arrow["\iota", shorten <=5pt, shorten >=5pt, Rightarrow, from=0, to=1-2]
		\end{tikzcd}\]
		\item (associativity)
		\[\begin{tikzcd}
		{G\times G\times\mathcal{X}} & {G\times\mathcal{X}} \\
		{G\times\mathcal{X}} & {\mathcal{X}}
		\arrow["{m\times \id_{\mathcal{X}}}", from=1-1, to=1-2]
		\arrow["{\id_G\times\mu}"', from=1-1, to=2-1]
		\arrow["\mu", from=1-2, to=2-2]
		\arrow["\alpha", shorten <=18pt, shorten >=18pt, Rightarrow, from=2-1, to=1-2]
		\arrow["\mu", from=2-1, to=2-2]
		\end{tikzcd}\]
	\end{enumerate}
    We call $\iota$ the \emph{unitor} and $\alpha$ the \emph{associator}.
    For $x\in\cX(U)$ and $g,h,k\in G(U)$, the components $\iota_x:x\to e\cdot x$ and $\alpha_{g,h,x}:g\cdot(h\cdot x)\to (gh)\cdot x$ satisfies the pentagon
    \[
    \alpha_{gh,k,x}\circ\alpha_{g,h,k\cdot x}=\alpha_{g,hk,x}\circ(g\cdot\alpha_{h,k,x}),
    \]
    and the unit identities are
    \[
    \alpha_{e,g,x}\circ\iota_{g\cdot x}=\id_{g\cdot x},\qquad\alpha_{g,e,x}\circ(g\cdot\iota_x)=\id_{g\cdot x}.
    \]
	We call $(\mathcal{X},\mu,\iota,\alpha)$ a \emph{$G$-prestack}.
    The action is \emph{strict} if both $\iota$ and $\alpha$ are identities.
\end{definition}
\begin{remark}
    Ginot and Noohi give the weak-action axioms for fibred groupoids and develop the quotient construction for topological stacks \cite[Definitions~3.1--3.3 and~4.2, \S 4.3]{ginotnoohi}. Here the acting group is the sheaf represented by the Lie group $G$ on the smooth site. Viewed as a groupoid over a point, also fits the framework of stacky groupoid actions \cite[\S 3]{bursztyn2020principal}. See \cite[Definitions~3.4--3.5]{zhu2009n} for the general stacky Lie groupoid convention.
\end{remark}
\begin{remark}\label{rem:pseudofunctor}
    For $U$ a manifold, the action is by $G(U)=C^\infty(U,G)$ on $\cX(U)$, naturally in $U$. The categorical description starts with the one-object groupoid of a group \cite[Example~1.3.9]{riehl2017category}. For each manifold $U$, let $\mathsf{B}(G(U))$ denote the one-object groupoid with automorphism group $G(U)$, regarded as a strict $2$-category with only identity $2$-morphisms. A $G$-action on $\mathcal X$ amounts to a family of pseudofunctors
    \[
        \mathfrak X_U:\mathsf B(G(U))\longrightarrow\mathbf{Gpd},
        \qquad
        \mathfrak X_U(*)=\mathcal X(U),
        \qquad
        \mathfrak X_U(g)=\mu_{U,g}:\mathcal X(U)\to\mathcal X(U),
    \]
    compatible with base change, i.e., for every map $f:V\to U$ and every $g\in G(U)$, one specifies a natural isomorphism $f^*\circ\mu_{U,g}\Rightarrow\mu_{V,g\circ f}\circ f^*$.
    These isomorphisms satisfy the coherence conditions for
    identities and composition of base maps, using the
    pseudofunctoriality constraints of $\mathcal X$, and are
    compatible with the action constraints.

    The same data can be encoded by a single pseudofunctor on an ordinary category. Let $\mathcal B_G$ be the Grothendieck construction of the prestack $U\mapsto\mathsf B(G(U))$. Its objects are manifolds, and an arrow $(f,g):V\to U$ consists of a smooth map $f:V\to U$ and a section $g\in G(V)$. For $(h,k):W\to V$, composition and identities are given by
    \[
    (f,g)\circ(h,k)
    =
    \bigl(f\circ h,(g\circ h)k\bigr),
    \qquad
    \operatorname{id}_U=(\operatorname{id}_U,e_U).
    \]
    The action determines a pseudofunctor
    \[
    \mathfrak X:\mathcal B_G^{\mathrm{op}}\longrightarrow\mathbf{Gpd},
    \qquad
    \mathfrak{X}(U)=\mathcal X(U), \quad \mathfrak X(f,g)=\mu_{V,g^{-1}}\circ f^*
    :\mathcal X(U)\longrightarrow\mathcal X(V).
    \]
    The inverse compensates for contravariance, and the coherence isomorphisms of $\mathfrak X$ are induced by those of $\mathcal X$ and the action. Let $\jmath:\mathbf{Man}\to\mathcal B_G$ be the functor which is the identity on objects and sends $f:V\to U$ to $(f,e_V)$. The unit constraint gives a pseudonatural equivalence
    \[
    \mathfrak X\circ \jmath^{\mathrm{op}}\simeq\mathcal X.
    \]
    Conversely, a pseudofunctor $\mathfrak X$ together with such an identification recovers the action: after transporting along this identification, $\mu_{U,g}$ is given by $\mathfrak X(\operatorname{id}_U,g^{-1})$.
    This is the smooth analogue of the formulation in \cite[Remark 1.4]{romagny2005group}, where the corresponding pseudofunctor is called a lax presheaf.
\end{remark}

\begin{definition}
A \emph{$G$-equivariant morphism} $(f,\sigma):(\cX,\mu,\iota,\alpha)\to(\cX',\mu',\iota',\alpha')$ consists of a morphism $f:\cX\to\cX'$ and an invertible $2$-cell
	\[\begin{tikzcd}
	{G\times\mathcal{X}} & {G\times\mathcal{X}'} \\
	{\mathcal{X}} & {\mathcal{X}'}
	\arrow["{\id_G\times f}", from=1-1, to=1-2]
	\arrow["\mu"', from=1-1, to=2-1]
	\arrow["\sigma", shorten <=25pt, shorten >=18pt, Rightarrow, from=2-1, to=1-2]
	\arrow["{\mu'}", from=1-2, to=2-2]
	\arrow["f", from=2-1, to=2-2]
	\end{tikzcd}\]
    For $x\in\cX(U)$ and $g,h\in G(U)$, the components $\sigma_{g,x}:f(g\cdot x)\to g\cdot f(x)$ satisfy
    \[
    \sigma_{gh,x}\circ f(\alpha_{g,h,x})
    =\alpha'_{g,h,f(x)}\circ(g\cdot\sigma_{h,x})\circ\sigma_{g,h\cdot x},
    \qquad
    \sigma_{e,x}\circ f(\iota_x)=\iota'_{f(x)}.
    \]
    A $2$-morphism $\eta:(f,\sigma)\Rightarrow (f',\sigma')$ is \emph{$G$-equivariant} when $(g\cdot\eta_x)\circ\sigma_{g,x}=\sigma'_{g,x}\circ\eta_{g\cdot x}$.
\end{definition}

The same definitions apply to stacks. We write $\mathbf{PSt}_G$ and $\mathbf{St}_G$ the $(2,1)$-categories of $G$-prestacks and $G$-stacks, respectively.

\begin{definition}\label{def:triv-forget}
    The trivial-action and forgetful functors are
    \[
    \begin{aligned}
    \Trivial:\St&\longrightarrow\St_G,&
    \cX&\longmapsto(\cX,\operatorname{pr}_2,\id,\id),\\
    \Forgetful:\St_G&\longrightarrow\St,&
    (\cX,\mu,\iota,\alpha)&\longmapsto\cX.
    \end{aligned}
    \]
    On $1$-morphisms $\Trivial(f)=(f,\id)$ and $\Forgetful(f,\sigma)=f$. Both send a $2$-morphism to its underlying $2$-morphism. We use the same notation for prestacks.
\end{definition}

\begin{remark}\label{rem:triv-forget}
There is a canonical identification $\Forgetful\Trivial\simeq\id_{\St}$. 
\end{remark}

\begin{proposition}[Adjoints of the forgetful functor]\label{prop:forgetful-adjoints}
In the $(2,1)$-category of stacks, the forgetful functor has a left and a right biadjoint. They are the free-action construction $\cY\mapsto G\times\cY$, with left translation on $G$, and the coinduced-action construction $\cY\mapsto\underline{\Hom}(G,\cY)$, with action
\[
(g\cdot f)(a)=f(ag).
\]
Here $\underline{\Hom}(G,\cY)(U)=\Mor_{\St}(G\times U,\cY)$, and the formula for the action is interpreted on smooth families. More precisely, $G\times(-)\dashv\Forgetful\dashv\underline{\Hom}(G,-)$, i.e.
\[
\begin{aligned}
\Mor_{\St_G}(G\times\cY,\cX)
&\simeq\Mor_{\St}(\cY,\Forgetful(\cX)),\\
\Mor_{\St_G}(\cX,\underline{\Hom}(G,\cY))
&\simeq\Mor_{\St}(\Forgetful(\cX),\cY).
\end{aligned}
\]
\end{proposition}
\begin{proof}
    For the first equivalence, given an equivariant morphism $(f,\sigma):G\times\cY\to\cX$ we restrict it to $\{e\}\times\cY$, and obtain $f|_{\{e\}\times\cY}:\cY\to\Forgetful(X)$. Its quasi-inverse sends $f:\cY\to\Forgetful(\cX)$ to $(g,y)\mapsto g\cdot f(y)$, with equivariance supplied by the action associator of $\cX$.
    
    For the second equivalence, take $f:\cX\to\underline{\Hom}(G,\cY)$, with $f(x):G\times U\to\cY$ for $x\in\cX(U)$, then evaluating at $e\in G$, we obtain a map $f(-)(e):\Forgetful(\cX)\to\cY$, where $f(x)(e):U\to\cY$ defines an element of $\cY(U)$. Its inverse sends $f:\Forgetful(\cX)\to\cY$ to the morphism which assigns to $x:U\to\cX$ the morphism defined by $(a,u)\mapsto f(a\cdot x)(u)$. The associator gives its equivariance, since evaluation at $a$ of the translate by $g$ is evaluation of the original family at $ag$.
    
    These constructions act on $2$-morphisms by restriction and by applying the corresponding formulas, and commute with precomposition and postcomposition up to the canonical coherence. This gives the two biadjunctions.
\end{proof}

\begin{definition}[{\cite[Definition 2.9]{barbosa2021equivariant}}]
    A \emph{differentiable $G$-stack} is a differentiable stack with a smooth $G$-action. A \emph{$G$-atlas} is a representable $G$-equivariant morphism $p:X\to\cM$, with $X$ a smooth $G$-manifold, such that $p$ is a surjective submersion.
\end{definition}

We denote with $\mathbf{DiffSt}_G$ the $(2,1)$-category of differentiable $G$-stacks.

\begin{proposition}[Existence of equivariant presentations]\label{prop:equivariant-atlas}
    Let $\cM$ be a differentiable $G$-stack. If $p':U_0\to\cM$ is any atlas, then
    \[
    p:X_0=G\times U_0\longrightarrow\cM,\qquad (g,u)\longmapsto g\cdot p'(u)
    \]
    is a $G$-atlas, for the $G$-action $h\cdot(g,u)=(hg,u)$. Its induced groupoid $\rX$ admits an induced strict smooth $G$-action.
\end{proposition}
\begin{proof}
Let $\operatorname{str}(g,x)=(g,g\cdot x)$. This is an isomorphism of stacks with quasi-inverse $(g,y)\mapsto(g,g^{-1}\cdot y)$ and comparison cells supplied by $\alpha,\iota$. The factorization is
\[
  G\times U_0\xrightarrow{\id_G\times p'}G\times\cM
    \xrightarrow{\operatorname{str}}G\times\cM
    \xrightarrow{\operatorname{pr}_2}\cM .
\]
The outer maps are representable surjective submersions. Equivariance for $p$ is the cell $\alpha_{h,g,p'(u)}^{-1}:p(hg,u)\to h\cdot p(g,u)$. Writing this cell as $\lambda_{h,z}$, an atlas arrow $a:p(z)\to p(z')$ is sent to $\lambda_{h,z'}^{-1}(h\cdot a)\lambda_{h,z}$. Naturality of $\alpha$ and the equivariance axiom for $\lambda$ give the strict group law on these arrows. The passage from a given equivariant atlas to the action on its nerve is \cite[Proposition~2.11]{barbosa2021equivariant}.
\end{proof}
\begin{remark}
    This fact is already noted in \cite[\S 3.1]{lerman2012hamiltonian}. Romagny proves a categorical strictification result for fibred groupoids in \cite[Proposition~1.5]{romagny2005group}.
\end{remark}

For a $G$-stack $\cX$ define $[\cX/G]$ by stackifying the transformation prestack. Over $T$ its objects are $x\in\cX(T)$ and an arrow $x\to y$ is $(g,a)$ with $g:T\to G$ and $a:g\cdot x\to y$. The unit is $(e,\iota_x^{-1})$, and the composite of $(g,a)$ followed by $(h,b)$ is
\[
    (hg,\ b\circ(h\cdot a)\circ\alpha_{h,g,x}^{-1}).
\]
The map $q:\cX\to[\cX/G]$ uses degree-$e$ arrows. The degree defines $r:[\cX/G]\to\BG$. Equivalently a quotient object is a principal right $G$-bundle with an equivariant morphism to $\cX$, its right action converted to a left action by inversion. See \cite[Definition~2.3(ii) and Proposition~2.6]{romagny2005group} and \cite[\S\S 4.1--4.2]{ginotnoohi}.

\begin{proposition}
    If $\cM$ is a differentiable $G$-stack, then $[\cM/G]$ is differentiable and $q:\cM\to[\cM/G]$ is a principal $G$-bundle.
\end{proposition}
\begin{proof}
    Apply \cite[Corollary~5.3]{bursztyn2020principal} to the one-object Lie groupoid $\mathsf{B}G=(G\rightrightarrows\pt)$.
\end{proof}

\begin{example}[Actions on a classifying stack]\label{ex:classifying-extension}
    Fix a Lie group $K$ and the trivial right principal $K$-bundle $P\cong K\to \pt$. A $G$-action on $\mathcal B K$ gives an extension
    \[
    \begin{tikzcd}
    1 & K & \widetilde G & G & 1
    \arrow[from=1-1,to=1-2]
    \arrow[from=1-2,to=1-3]
    \arrow["\pi",from=1-3,to=1-4]
    \arrow[from=1-4,to=1-5]
    \end{tikzcd}
    \]
    where $\widetilde G=\Set{(g,\varphi)|g\in G, \ \varphi:g\cdot P\to P\text{ isomorphism of principal $K$-bundle}}$ with multiplication defined as
    \[
    (g,\varphi)(h,\psi)=(gh,\varphi\circ(g\cdot\psi)\circ\alpha_{g,h,P}^{-1}).
    \]
    This is the isotropy Lie group of $q(P)$ in $[\mathcal B K/G]$. Its kernel identifies with $K$ by left translations after the identification at $e$ supplied by $\iota_P$. To check surjectivity in families, take $g:T\to G$. The torsor $g\cdot P_T$ on $T$ is locally trivial, so a trivialization supplies a smooth isomorphism $g\cdot P_T\to P_T$. Thus $\pi$ is an epimorphism of sheaves. Taking $T=G$ and $g=\id_G$ produces smooth local sections, and translation in the Lie groups makes $\pi$ a submersion at every point.

    Conversely, an extension as above acts on $G$ by $\widetilde g\ast a=a\pi(\widetilde g)^{-1}$. The action is transitive with stabilizer $K$ at $e$, so $[G/\widetilde G]\simeq\mathcal B K$. Left translation by $G$ commutes with this action and therefore induces an action on this quotient stack; finally, transport across this equivalence gives a $G$-action on $\mathcal B K$.
    
    The map $\pt\to\mathcal B K$ defined by $P$, with trivial $G$-action on $\pt$, is $G$-equivariant precisely when there is a smooth family $\beta_g:g\cdot P\to P$ satisfying the action cocycle. Such a family is exactly a smooth homomorphic section $g\mapsto(g,\beta_g)$ of $\pi$. Thus this atlas is equivariant if and only if the extension splits.
    See \cite[Example 3.4]{lerman2012hamiltonian}.
\end{example}

Let us review some classes of differentiable stacks that we will be using in the course of this discussion.
\begin{definition}
    Let $\cM\simeq[X_0/X_1]$ be a differentiable stack.
    \begin{itemize}
        \item It is \emph{separated} if it admits a proper Lie groupoid presentation, meaning that $(s,t):X_1\to X_0\times X_0$ is proper, or equivalently, the diagonal is proper in the quasi-representable sense. Properness is Morita invariant.
        \item It is \emph{étale} (or \emph{Deligne--Mumford}) if it admits an étale presentation, i.e. one with $s,t$ local diffeomorphisms;
        \item It is \emph{regular} if the isotropy dimensions are locally constant. In a presentation this is equivalent to locally constant rank of the algebroid anchor $\rho:A\to TX_0$. Equivalently, the orbit dimensions are locally constant on each component of $X_0$.
        \item It is an \emph{orbifold stack} if it admits a proper étale presentation.
    \end{itemize}
\end{definition}
\begin{remark}
    Moerdijk describes a regular Lie groupoid as an extension of a foliation groupoid by the bundle of identity components of its isotropy groups \cite[\S~2]{moerdijk2002classification}.
\end{remark}

\begin{proposition}\label{prop:proper-quotient}
    Let $\cM$ be a differentiable $G$-stack. If $G$ is compact and $\cM$ is separated, then $[\cM/G]$ is separated.
\end{proposition}
\begin{proof}
    Choose the $G$-atlas of \Cref{prop:equivariant-atlas}. Its Lie groupoid $\rX$ is proper by Morita invariance. The composite $X_0\to[\cM/G]$ is an atlas, with arrows $(g,\gamma)$ interpreted as $\gamma:g\cdot x\to y$. Thus the induced groupoid $\rY$ has
    \[
    Y_1=G\times X_1,\qquad s_{\rY}(g,\gamma)=g^{-1}\cdot s(\gamma),\qquad t_{\rY}(g,\gamma)=t(\gamma).
    \]
    The degree is $\rr_1(g,\gamma)=g$ and composition is $(h,\delta)\circ(g,\gamma)=(hg,\delta\circ(h\cdot\gamma))$. The map
    \[
    G\times X_1\longrightarrow G\times X_0\times X_0,
    \quad(g,\gamma)\longmapsto(g,g^{-1}\cdot s(\gamma),t(\gamma))
    \]
    is the product of the proper anchor of $\rX$ with $G$, followed by a diffeomorphism of its target, and is therefore proper. Projection to $X_0\times X_0$ is proper because $G$ is compact. Their composite is the anchor of $\rY$.
\end{proof}

\begin{lemma}[A slice compatible with the map to $\BG$]\label{lem:relative-slice}
    Let $\cN$ be a separated differentiable stack, $r:\cN\to\BG$ be a morphism, and let $y$ be a point with compact isotropy $L=\Aut_{\cN}(y)$. Choose a trivialization of the principal $G$-bundle $r(y)$. Then there are an open neighbourhood $\mathcal{U}\subseteq\cN$ of $y$, an $L$-invariant open ball $V$ about $0$ in the normal representation at $y$, and an isomorphism
    \[
    \mathcal{U}\simeq[V/L]
    \]
    under which $r|_{\mathcal{U}}$ is $2$-isomorphic to the morphism induced by the isotropy homomorphism $\rho:L\to G$. If $\rho$ is surjective and $K=\ker\rho$, then
    \[
    \mathcal{U}\times_{\BG}\pt\simeq[V/K],
    \]
    with the residual $G$-action induced by 
    \begin{equation}\label{eq:slice-extension}
        \begin{tikzcd}
        	1 & K & L & G & 1.
        	\arrow[from=1-1, to=1-2]
        	\arrow[from=1-2, to=1-3]
        	\arrow["\rho", from=1-3, to=1-4]
        	\arrow[from=1-4, to=1-5]
        \end{tikzcd}
    \end{equation}
\end{lemma}
\begin{proof}
    Take a proper presentation of $\cN$ and a slice at a representative of $y$. By \cite[Lemma 3.10]{crainic2013linearization}, restriction to the slice is Morita equivalent to restriction to its open saturation. Linearization at the fixed point of the slice groupoid, followed by shrinking to an invariant ball, gives $\mathcal{U}\simeq[V/L]$. The groupoid results used here are \cite[Theorem 2.1, Lemma 3.10 and Proposition 3.11]{crainic2013linearization}; see also the explicit construction in \cite[Section 5.3]{Mrcun2018} and \cite[Theorem 3.7.2]{cnossen2024twisted} for a precise statement on the local structure of separated differentiable stacks.

    The pullback of $r$ to $[V/L]$ is a $L$-equivariant principal $G$-bundle $P\to V$, and we have that the left $L$-action commutes with the right $G$-action. So, if we choose any $p_0\in P_0$ in the fiber of $0\in V$, the equation $l p_0=p_0\rho(l)$ defines the given homomorphism $\rho$. Choose a principal connection on $P$ and average it.
    The affine conditions defining a principal connection are preserved by this average, which is $L$-invariant. Parallel transport along $t\mapsto tv$, $0\leq t\leq1$, gives a smooth trivialization
    \[
    P\cong V\times G,\qquad
    l\cdot(v,a)=(l\cdot v,\rho(l)a).
    \]
    Indeed, the action sends the radial path to the radial path for $l\cdot v$, and preserves both the connection and right translation. This proves compatibility of the slice with $r$.

    The $2$-fibre product with $\pt\to\BG$ is presented by $[P/L]=[(V\times G)/L]$. If $\rho$ is surjective, $V\times\{e\}$ meets all $L$-orbits and its restricted arrows are exactly $K\times V$. The map
    \[
    L\times V\longrightarrow V\times G,
    \qquad(l,v)\longmapsto(l v,\rho(l))
    \]
    is a surjective submersion, so this restriction is a Morita equivalence. Right translation on the $G$-coordinate, written as a left action using inverses, is the residual $G$-action; descent identifies it with the action associated to the extension, and so these constructions supply the asserted equivalence together with its $2$-isomorphism over $\BG$. See also \cite[Example 4.1]{lerman2012hamiltonian} for further details.
\end{proof}

We will need the compact-open topology on the space of homomorphisms and a criterion for its smooth test families. The classical conjugacy and deformation results include \cite{montgomery1940topological,conner1962differentiable,weil1964remarks,nijenhuis1967deformations,lee1969conjugacy,lee1974deformations}.

\begin{theorem}\label{thm:homspace}
    Let $G$ be a compact Lie group and $K$ any Lie group. Give the space of smooth homomorphisms $\operatorname{hom}(G,K)$ the compact-open topology. For each homomorphism $\varphi$ set
    \[
    C_K(\varphi)=\Set{k\in K | k\varphi(g)=\varphi(g)k\ \forall\,g\in G}.
    \]
    Then every conjugacy orbit in $\operatorname{hom}(G,K)$ is clopen, and conjugation induces a homeomorphism
    \[
    \coprod_{[\varphi]\in\operatorname{hom}(G,K)/K}K/C_K(\varphi)
    \;\xrightarrow{\ \cong\ }\;\operatorname{hom}(G,K).
    \]
    The set of conjugacy classes is countable, so, via this homeomorphism, the disjoint union on the left gives $\operatorname{hom}(G,K)$ a smooth manifold structure.
    
    For every manifold $T$, a map $a:T\to\operatorname{hom}(G,K)$ is smooth for this structure if and only if its evaluation map
    \[
    \operatorname{ev}_a:T\times G\longrightarrow K, \qquad (t,g)\longmapsto a(t)(g)
    \]
    is smooth. In particular, any smooth family $\varphi_t$ is locally conjugate to a fixed homomorphism by a smooth family, i.e. for every $t_0\in T$ there exist a neighborhood $U\ni t_0$ and a smooth map $k:U\to K$, with $k(t_0)=e$, such that
    \[
    \varphi_t(g)=k(t)\varphi_{t_0}(g)k(t)^{-1}.
    \]
\end{theorem}
\begin{proof}
    Schwede proves that the conjugacy orbits of $K^0$ are open for compact source and arbitrary Lie target \cite[Proposition~A.25]{schwede}. The $K$-orbits are therefore clopen. Rezk gives the displayed homeomorphism under the same compactness hypothesis \cite[Proposition~7.1]{rezk2018classifying}. Since $G$ is compact and metrizable and $K$ is second countable and metrizable, $C(G,K)$ with the compact-open topology is second countable. One may take finite intersections of conditions mapping the closure of a basic open set of $G$ into a basic open set of $K$, using a countable relatively compact basis in $G$. The subspace of homomorphisms thus has at most countably many disjoint open orbits.
    
    Suppose that the evaluation map of $a$ is smooth. Choose coordinates near $t_0$ with $t_0=0$, and put $\varphi_{s,t}(g)=a(st)(g)$ for $0\leq s\leq1$ on a coordinate ball. Its right logarithmic derivative $A_{s,t}(g)=\partial_s\varphi_{s,t}(g)\,\varphi_{s,t}(g)^{-1}$ satisfies
    \[
     A_{s,t}(gh)=A_{s,t}(g)+\operatorname{Ad}_{\varphi_{s,t}(g)}A_{s,t}(h).
    \]
    For normalized Haar measure put $u_{s,t}=\int_G A_{s,t}(h)\,dh$. It depends smoothly on $(s,t)$ and $A_{s,t}(g)=u_{s,t}-\operatorname{Ad}_{\varphi_{s,t}(g)}u_{s,t}$. Solve $\partial_s k_{s,t}\,k_{s,t}^{-1}=u_{s,t}$ with $k_{0,t}=e$. After shrinking the ball, smooth dependence of the solution along the constant solution at $t=0$ gives $k_{s,t}$ for $0\leq s\leq1$. Differentiation makes $k_{s,t}^{-1}\varphi_{s,t}(g)k_{s,t}$ independent of $s$. Thus $k_{1,t}$ is the required smooth conjugator. This is the parameter version of the averaging argument of \cite[Theorems~1.3 and~3.8]{cardenas2020stability}. Conversely, a smooth map into the homogeneous orbit lifts locally through the submersion $K\to K/C_K(\varphi)$, and conjugation then gives a jointly smooth evaluation.
\end{proof}

\section{Fixed locus stack}\label{sec:fixed-locus}

Let $U$ be a manifold. If $\cX$ is a stack over $\Man$, its restriction to the slice site $\Man/U$ is the stack
\[
\cX_U:=\cX|_{\Man/U},
\qquad
\cX_U(V\xrightarrow{u}U):=\cX(V).
\]
If $x\in\cX(U)$, then $x$ determines a section of $\cX_U$: over $u:V\to U$ it is the pullback
\[
u^*x\in\cX(V).
\]
Thus, whenever we write \(x|_{(V\xrightarrow{u}U)},\) we mean $u^*x$.  For a smooth section $g\in G(V)$, the object \(g\acts x|_{(V\xrightarrow{u}U)}\) means \(g\acts u^*x=\mu_{g,V}(u^*x)\in\cX(V).\)

Let $G_U$ be the represented group sheaf on $\Man/U$:
\[
G_U(V\xrightarrow{u}U)=C^\infty(V,G).
\]

\subsection{Definition and properties}

Fix $x\in\cX(U)$. For $u:V\to U$ and $g\in G_U(V)$, put
\[
a_x(g)=g\cdot u^*x,\qquad b_x(g)=u^*x.
\]

\begin{definition}\label{def:stabilizer}
    The \emph{extended stabilizer sheaf} of $x$ is the group sheaf
    \[
    \Stab^G_{\cX}(x)(V\xrightarrow{u}U) =\Set{(g,\varphi) | g\in C^\infty(V,G),\ \varphi:g\cdot u^*x\xrightarrow{\sim}u^*x}.
    \]
    Equivalently, it is a sheaf representative of the following $2$-fibre product on the slice site:
    \[
    \begin{tikzcd}
    {\Stab^G_{\cX}(x)} & {\cX_U} \\
    {G_U} & {\cX_U\times\cX_U}.
    \arrow[from=1-1,to=1-2]
    \arrow["\alpha_x"',from=1-1,to=2-1]
    \arrow["\Delta_U",from=1-2,to=2-2]
    \arrow["{(a_x,b_x)}",from=2-1,to=2-2]
    \arrow["\lrcorner"{anchor=center,pos=0.125},draw=none,from=1-1,to=2-2]
    \end{tikzcd}
    \]
    The projection $\alpha_x(g,\varphi)=g$ is a homomorphism. Writing $z=u^*x$, its multiplication and unit are
    \[
     (g,\varphi)(h,\psi)
       =(gh,\varphi\circ(g\cdot\psi)\circ\alpha_{g,h,z}^{-1}),
     \qquad 1=(e,\iota_z^{-1}).
    \]
    Its kernel is the automorphism sheaf of $x$, with $a\mapsto(e,a\circ\iota_z^{-1})$ on sections.
\end{definition}
\begin{remark}
    Fixing arrows glue uniquely because the isomorphisms in $\cX$ form sheaves. The pentagon and unit axioms for the action give the group identities in $\Stab_{\cX}^G(x)$.
\end{remark}

The \emph{universal stabilizer}, see \cite[\S~4.1]{romagny2022algebraicity}, is the stack $\Stab^G_{\cX}$ defined by the $2$-Cartesian square
\[
\begin{tikzcd}
{\Stab^G_{\cX}} & {\cX} \\
{G\times\cX} & {\cX\times\cX}.
\arrow[from=1-1,to=1-2]
\arrow[from=1-1,to=2-1]
\arrow["\Delta",from=1-2,to=2-2]
\arrow["{(\mu,\operatorname{pr}_2)}",from=2-1,to=2-2]
\arrow["\lrcorner"{anchor=center,pos=0.125},draw=none,from=1-1,to=2-2]
\end{tikzcd}
\]
Its pullback along an object $x:U\to\cX$ is the group sheaf just defined.

\begin{definition}\label{def:fixed-locus}
    The \emph{fixed locus stack} $\cX^G\subseteq\cX$ is the full substack whose objects $x\in\cX(U)$ satisfy the condition that $\alpha_x:\Stab^G_{\cX}(x)\to G_U$ is an epimorphism of sheaves on $\Man/U$. Explicitly, for every $u:V\to U$ and every $g:V\to G$, there is an open cover $\{V_i\to V\}_{i\in I}$ with isomorphisms
    \[
    g|_{V_i}\cdot(u|_{V_i})^*x\xrightarrow{\sim}(u|_{V_i})^*x.
    \]
    These local choices need not be multiplicative in $g$.
\end{definition}

For prestacks the same condition defines the fixed locus prestack.

\begin{remark}\label{rem:algebraic-fixed-locus}
    This is the smooth, sheaf-theoretic analogue of \cite[Definition~A.9 and Remark~A.10]{aranha2025virtual}: the defining condition there is that the stabilizer homomorphism be an effective epimorphism of fppf sheaves. Here the topology is the open-cover topology on manifolds and an epimorphism of sheaves means local lifts of each section.
\end{remark}

\begin{proposition}\label{prop:fixed-formal}
    Let $\cX$ be a $G$-stack.
    \begin{enumerate}[label=(\alph*)]
        \item $\cX^G$ is a full substack of $\cX$.
        \item Being fixed is local on the test manifold.
        \item It is invariant under isomorphisms of objects.
        \item Every $G$-morphism $f:\cX\to\cY$ restricts to $f|_{\cX^G}:\cX^G\to\cY^G$.
        \item For $u:V\to U$, restriction gives
        \[
        \Stab^G_{\cX}(x)|_{\Man/V}\cong\Stab^G_{\cX}(u^*x),
        \]
        with the canonical pullback coherence understood.
    \end{enumerate}
\end{proposition}
\begin{proof}
    For (b), suppose $x$ is fixed on an open cover $U_i$. Given $g:V\to G$ over $u:V\to U$, first restrict to $u^{-1}(U_i)$ and then use local lifting there. Refining these covers supplies local lifts on $V$. Conversely the condition is preserved by restriction. An isomorphism $\eta:x\to y$ identifies stabilizers by
    \[
    (g,\varphi)\longmapsto(g,\eta\circ\varphi\circ(g\cdot\eta)^{-1}),
    \]
    which proves (c). Fullness is the definition, and effective descent for objects in $\cX$, together with (b) and (c), proves (a). For (d), apply $f$ to each local lift and use the equivariance $2$-cell to obtain a lift for $f(x)$. Finally pullback of an isomorphism and the action's naturality give (e). The same arguments respect $2$-morphisms, so equivariant equivalences induce equivalences on fixed substacks.
\end{proof}

\begin{remark}
    The inclusion is fully faithful. Even when its pullback to an atlas is an embedded submanifold, its pullback to an arbitrary test manifold may be singular. This distinction matters for quasi-representability and representability.
\end{remark}

\begin{proposition}[Functoriality and the induced action]\label{prop:fixed-categorical}
    The assignments
    \[
    \cX\longmapsto\cX^G,\quad (f,\sigma)\longmapsto f|_{\cX^G},\quad \eta\longmapsto\eta|_{\cX^G}
    \]
    define a pseudofunctor $\FixedLocus:\St_G\to\St$, with a canonical pseudonatural transformation
    \[
    i:\FixedLocus\Longrightarrow\Forgetful,
    \qquad i_{\cX}:\cX^G\hookrightarrow\cX.
    \]
    The fixed locus also carries the restricted $G$-action, for which $i_{\cX}$ is equivariant.
\end{proposition}
\begin{proof}
The restriction of morphisms and $2$-morphisms follows from \Cref{prop:fixed-formal}. The identity and composition constraints restrict because $\cX^G$ is a full substack, and this also gives pseudonaturality of $i$. If $x\in\cX^G(U)$ and $h\in G(U)$, the defining epimorphism implies that $h\cdot x$ is locally isomorphic to $x$. Isomorphism invariance and locality therefore imply that $h\cdot x$ belongs to $\cX^G(U)$. This restricts the action morphism; its coherence cells restrict by fullness. The extended stabilizer of a fixed object is unchanged by passage to the full invariant substack. Hence every object of the restricted-action stack is fixed.
\end{proof}

\begin{proposition}\label{prop:fixed-products}
    The functor $\FixedLocus$ preserves the terminal object and finite products, up to canonical equivalence. It does not in general preserve $2$-fiber products.
\end{proposition}
\begin{proof}
    Let $\cX,\cY$ be $G$-stacks. For the diagonal action, the extended stabilizer of $(x,y)$ is the fiber product over $G_U$ of the extended stabilizers of $x$ and $y$. If both stabilizer projections are epimorphisms, a common refinement of the two covers gives simultaneous local lifts. Conversely, a lift to the fiber product projects to a lift in each factor. By fullness this identifies the fixed groupoid of the product with the product of the fixed groupoids. The unique object of the terminal stack is fixed.

    For the second statement, consider $G=K=\mathbb Z/2$. Give $\pt$ and $\mathcal BK$ the trivial $G$-actions, and  let $f:\pt\to\mathcal{B}K$ classify the trivial principal $K$-bundle $x$ over $\pt$. Since $\Aut_{\mathcal B K}(x)\cong K$, a $G$-equivariant structure on the map $f$ (with both source and target actions trivial) is specified by a family of automorphisms $\sigma_g:x\xrightarrow{\sim}x$ satisfying 
    \[\sigma_{gh}=\sigma_h\circ\sigma_g, \qquad \sigma_0=\mathrm{id}_x.\]
    In particular $g\mapsto\sigma_g$ is a homomorphism here because $K=\mathbb Z/2$ is abelian.
    Choose the two homomorphisms
    \[
      \sigma_{triv}(g)=0,
      \qquad
      \sigma_{id}(g)=g.
    \]
    They define two equivariant morphisms $(f,\sigma_{triv}),(f,\sigma_{id}):\pt\to\mathcal BK$ whose underlying morphisms of stacks are identical, but their equivariance constraints are not isomorphic: in fact, an equivariant $2$-isomorphism between them would amount to $k\in K$ satisfying $\sigma_{id}(g)k=k\sigma_{triv}(g)$ for all $g\in G$, which is impossible.
    
    On the underlying $2$-pullback $K\simeq\pt\times_{\mathcal B K}\pt$, the equivariance convention is $(\sigma_i)_{g,x_i}:f_i(g\cdot x_i)\to g\cdot f_i(x_i)$. A gluing arrow $\varphi:f_0(x_0)\to f_1(x_1)$ is transported to
    \[
        \varphi_g=(\sigma_1)_{g,x_1}^{-1}\circ(g\cdot\varphi)\circ(\sigma_0)_{g,x_0}.
    \]
    Thus $g\cdot (x_0,x_1,\varphi)=(g\cdot x_0,g\cdot x_1,\varphi_g)$. In additive notation for $K=\mathbb Z/2$, the resulting action is $g\cdot k=-\sigma_{\rm id}(g)+k+\sigma_{\rm triv}(g)=k+g$.
    Thus the nonidentity element of $G$ exchanges the two points of $K$. This is a free action on the $2$-fiber product $K$, and therefore
    \[
      K \simeq \pt^{G}\times_{\mathcal{B}K^{G}}\pt^{G}\not\simeq K^{G}=\varnothing.
    \]
\end{proof}

\subsection{Description by a Lie groupoid and differentiability}

Choose a $G$-atlas $p:X_0\to\cM$ as in \Cref{prop:equivariant-atlas}, and let $\rX=(X_1\rightrightarrows X_0)$ be its strict smooth equivariant groupoid. Write $\mathcal O_x=\rt(\rs^{-1}(x))$ and $H_x=\Aut_{\rX}(x)$.

\begin{definition}
    The \emph{extended $G$-stabilizer} of $x\in X_0$ is
    \[
    \widehat G_x=\Set{(g,\gamma)|g\in G,\ \gamma:g\cdot x\to x\text{ in } X_1},
    \qquad\pi_x(g,\gamma)=g.
    \]
    It is the isotropy Lie group of $x$ in the quotient groupoid $\rY$ of \Cref{prop:proper-quotient}. It fits into the diagram
    \[
    \begin{tikzcd}
    {\widehat G_x} & {X_1} \\
    G & {X_0\times X_0}
    \arrow[from=1-1,to=1-2]
    \arrow["\pi_x"',from=1-1,to=2-1]
    \arrow["{(s,t)}",from=1-2,to=2-2]
    \arrow["{g\mapsto(g\cdot x,x)}",from=2-1,to=2-2]
    \arrow["\lrcorner"{anchor=center,pos=0.125},draw=none,from=1-1,to=2-2]
    \end{tikzcd}
    \]
    with multiplication $(g,\gamma)(h,\delta)=(gh,\gamma\circ(g\cdot\delta))$.
\end{definition}

The kernel is $H_x$. Whenever $\pi_x$ is surjective, we obtain the extension
\begin{equation}\label{eq:extended-isotropy}
\begin{tikzcd}
1 & {H_x} & {\widehat G_x} & G & 1.
\arrow[from=1-1,to=1-2]
\arrow[from=1-2,to=1-3]
\arrow["\pi_x",from=1-3,to=1-4]
\arrow[from=1-4,to=1-5]
\end{tikzcd}
\end{equation}
With the chosen lift $x$, $\widehat G_x\cong\Aut_{[\cM/G]}(q(p(x)))$. It represents the stabilizer sheaf of the object $p(x)$ over a point. If $\cM$ is separated and $G$ compact, it is compact by \Cref{prop:proper-quotient}.

\begin{lemma}
\label{lem:fixed-locus-differentiable}
    The fixed-object set in this atlas is
    \[
     X_0^G=\{x\in X_0\mid p(x)\in\cM^G(\pt)\}=\{x\in X_0\mid\pi_x\text{ is surjective}\}=\{x\in X_0\mid G\cdot x\subseteq\mathcal O_x\}.
    \]
    It is saturated. Suppose in addition that $X_0^G\subseteq X_0$ is an embedded submanifold and that the object $p|_{X_0^G}\in\cM(X_0^G)$ satisfies the family lifting condition
    \begin{equation}\label{eq:family-lifting}
    \Stab^G_{\cM}(p|_{X_0^G})\longrightarrow G_{X_0^G}
    \quad\text{is an epimorphism of sheaves.}
    \end{equation}
    Then
    \[
    X_0^G\simeq X_0\times_{\cM}\cM^G,
    \qquad
    \cM^G\simeq[X_0^G/(X_1|_{X_0^G})],
    \]
    where $X_1|_{X_0^G}=\rs^{-1}(X_0^G)=\rt^{-1}(X_0^G)$ and $\rX^G=(X_1|_{X_0^G}\rightrightarrows X_0^G)$ denote the presenting groupoid. In particular $\cM^G$ is differentiable and the inclusion is a quasi-representable embedding; it is closed if $X_0^G$ is closed.
\end{lemma}
\begin{proof}
    Over a point, the existence of an arrow $g\cdot x\to x$ for every $g$ is exactly surjectivity of $\pi_x$. Conversely a surjective Lie group homomorphism is a submersion and has local smooth sections; pulling these back along every $V\to G$ gives the required epimorphism of sheaves over $\pt$. This last smooth-lifting step is needed with the represented group convention. The orbit criterion follows immediately. If $\eta:y\to x$, an arrow $g\cdot x\to x$ transports to
    \[
    g\cdot y\xrightarrow{g\cdot\eta}g\cdot x\longrightarrow x\xrightarrow{\eta^{-1}}y,
    \]
    which proves saturation.

    Since the inclusion is fully faithful, the sections of $X_0\times_{\cM}\cM^G$ over $U$ are the maps $u:U\to X_0$ with $u^*p\in\cM^G(U)$. Such a map has image in $X_0^G$, and factors smoothly through $X_0^G$ because $X_0^G$ is embedded. Conversely, a map through $X_0^G$ pulls back \eqref{eq:family-lifting}, and hence gives a fixed object.
    
    Saturation gives $\rs^{-1}(X_0^G)=\rt^{-1}(X_0^G)$, an embedded manifold. The restrictions of $\rs,\rt$ are surjective submersions by base change. Thus the full restriction is a Lie groupoid. The map $X_0^G\to\cM^G$ is the pullback of the atlas $X_0\to\cM$, so it is a representable surjective submersion and yields the asserted presentation. Pulling the inclusion back along $p:X_0\to\cM$ gives the embedding $X_0^G\hookrightarrow X_0$.
\end{proof}

The full restriction above keeps every arrow of $X_1$ between its objects. By comparison, the strict fixed subgroupoid of an action on a presentation keeps only objects and arrows fixed as equalities by every $g$. For example, the atlas $G\times U_0$ of \Cref{prop:equivariant-atlas} carries a free action on objects. For nontrivial $G$ its strict fixed subgroupoid is empty even when the original stack action is trivial. The full restriction instead presents the intrinsic fixed stack under the family hypothesis of \Cref{lem:fixed-locus-differentiable}.

\begin{lemma}\label{lem:saturation}
    If $\rX$ is proper, then $X_0^{G}$ is a closed saturated subset of $X_0$.
\end{lemma}
\begin{proof}
    Since the groupoid is proper, the orbit relation
    \[
    R=(\rs,\rt)(X_1)\subseteq X_0\times X_0
    \]
    is closed.  Consider
    \[
    A=\Set{(g,x)\in G\times X_0 | (g\acts x,x)\in R}.
    \]
    This is closed in $G\times X_0$.  Its complement is open, and the complement of the fixed object locus is the projection of $(G\times X_0)\setminus A$ to $X_0$.  Since projections from a product are open, this complement is open.  Hence $X_0^{G}$ is closed. Saturation follows from the orbit criterion in \Cref{lem:fixed-locus-differentiable}.
\end{proof}

The correct local model must retain the morphism $[\cM/G]\to\BG$. By \Cref{lem:relative-slice}, near $q(x)$ the quotient stack has a chart $[V_x/\widehat G_x]$ compatible with $\pi_x:\widehat G_x\to G$. For a point of the fixed locus, pulling this chart back to $\cM$ gives $[V_x/H_x]$, with residual $G$-action induced by \eqref{eq:extended-isotropy}. Here $V_x$ is an invariant ball in the normal representation, not an open subset of the original object manifold in general.

This pointwise criterion and saturation agree with \cite[Definition~4.7 and Lemma~4.8]{jiang2026equivariant}. In their Kan presentation the degree map $\rr_1:Y_1\to G$ has $(\rs_{\rY},\rr_1)$ a surjective submersion \cite[Definition~2.4]{jiang2026equivariant}. The following theorem identifies the fixed substack on smooth families.

\begin{theorem}[Differentiability of the fixed locus]\label{thm:fixed-locus-proper-regular-stack}
    Let $G$ be compact and connected and let $\cM$ be a separated differentiable $G$-stack. Write $|\cM|^G\subseteq|\cM|$ for the pointwise fixed set, with the subspace topology of the coarse space. Assume that
    \[
           |\cM|^G\longrightarrow\mathbb N,\qquad [x]\longmapsto\dim(\Aut_{\cM}(x))
    \]
    is locally constant, for example when $\cM$ is regular. Then $\cM^G$ is a separated regular differentiable stack and its inclusion $i:\cM^G\to\cM$ is a quasi-representable closed embedding. For any atlas the inverse image of $|\cM|^G$ is closed, embedded and saturated, represents the atlas pullback of $i$, and its full restriction presents $\cM^G$.
\end{theorem}
\begin{proof}
    Fix a point $x$ of the fixed locus and use the relative slice above. Put $L=\widehat G_x$, $K=H_x$, $V=V_x$, and $L^0$ for the identity component of $L$, thus we obtain the exact sequence \eqref{eq:slice-extension}, with $L$ compact, and the local stack is $[V/K]$. Since a surjective Lie group homomorphism is a submersion, $\pi_x(L^0)$ is an open connected subgroup of $G$; connectedness gives $\pi_x(L^0)=G$.

    For $v\in V$ write $L_v$ for its stabilizer and $K_v=K\cap L_v$. Pointwise fixedness is $\pi_x(L_v)=G$. The map $V\to|\cM|$ is continuous, so the hypothesis at $[x]$ allows us to shrink the ball until $\dim K_v=\dim K$ for all such $v$. For each of these $v$, the identity component $L_v^0$ also maps onto $G$, consequently
    \[
    \dim L_v^0=\dim G+\dim K_v=\dim G+\dim K=\dim L^0.
    \]
    The closed connected subgroup $L_v^0\subseteq L^0$ is therefore all of $L^0$. Conversely, if $L^0$ fixes $v$, the stabilizer $L_v$ maps onto $G$. Thus the pointwise fixed subset is exactly
    \[
    F_V:=\{v\in V\mid\pi_x(L_v)=G\}=V^{L^0}.
    \]
    This is a closed embedded linear fixed submanifold intersected with the ball $V$.

    To check the family condition, let $v:T\to F_V$ and $g:T\to G$ be smooth, and pull back local smooth sections of $L^0\to G$ along $g$; they lift $g$ locally to elements of $L^0$ which fix $v$, hence to the required arrows of the extended stabilizer. Conversely, every fixed family in the atlas has image in $F_V$ and factors smoothly through this embedded submanifold. Locally trivializing the $K$-torsor of a general family gives the fixed stack $[F_V/K]$ in the slice and the atlas pullback $F_V\hookrightarrow V$.

    To pass to any atlas $A\to\cM$, compare $V$ with $A$ over the open slice by $D=A\times_{\cM}V$. Its projections $a:D\to A$ and $b:D\to V$ are submersions and $a$ covers the part of $A$ in this slice, and the inverse image $b^{-1}(F_V)$ is embedded and is saturated for $a$ because the fixed condition is invariant under isomorphisms. Along this submanifold its tangent contains $\ker(da)$, and every local section of $a$ is therefore transverse to it. Its inverse image along such a section is exactly the pointwise fixed set on the section's domain and is embedded. The family condition descends through the same sections. At points outside $|\cM|^G$ there is an open complement by \Cref{lem:saturation}, first on a $G$-atlas and hence on the orbit space by the quotient topology. The full restrictions glue, since each represents the same condition on test families.
    
    The restricted groupoid anchor is proper by base change of the original proper anchor to the square of the closed fixed-object subset, and fullness preserves all isotropy groups and their dimensions are locally constant on that subset by hypothesis, proving regularity.
\end{proof}

\begin{corollary}[Orbifold fixed locus]\label{cor:fixed-locus-orbifold-stack}
    Let $G$ be a compact connected Lie group acting smoothly on an orbifold stack $\mathcal M$. Then $\mathcal M^G$ is an orbifold stack, and $i:\mathcal M^G\to\mathcal M$ is a quasi-representable closed embedding. For every proper \'etale presentation, the inverse-image fixed-object subset is closed, embedded and saturated, and the full restricted groupoid presents $\mathcal M^G$.
\end{corollary}
\begin{proof}
    For an \'etale presentation all isotropy Lie algebras vanish. Apply \Cref{thm:fixed-locus-proper-regular-stack}: restricting source and target to the full saturated subgroupoid preserves local diffeomorphisms, and the restricted anchor is proper. Pointwise embeddedness is precisely \cite[Theorem~4.9]{jiang2026equivariant}. The family lifting calculation above identifies their full subgroupoid as a presentation of the sheaf-theoretic fixed stack.
\end{proof}

\section{Homotopy fixed locus stack}\label{sec:homotopy-fixed-locus}

\subsection{Definition and properties}

The homotopy fixed locus remembers a \emph{smooth coherent choice} of fixing isomorphisms.

\begin{definition}\label{def:homotopy-fixed}
    For a smooth $G$-action on a stack $\cX$, the \emph{homotopy fixed locus stack} is characterized by the natural equivalence
    \[
    \Mor_{\St}(\cY,\cX^{hG})
    \simeq\Mor_{\St_G}(\Trivial(\cY),\cX).
    \]
    Explicitly, an object over $U$ is $(x,\beta)$, where $x\in\cX(U)$ and
    \[\begin{tikzcd}
    	{G\times U} & U \\
    	{G\times\cX} & \cX
    	\arrow["{\operatorname{pr}_U}", from=1-1, to=1-2]
    	\arrow["{\id_G\times x}"', from=1-1, to=2-1]
    	\arrow["x", from=1-2, to=2-2]
    	\arrow["\beta", between={0.4}{0.6}, Rightarrow, nfold, from=2-1, to=1-2]
    	\arrow["\mu", from=2-1, to=2-2]
    \end{tikzcd}\]
    is a 2-isomorphism over $G\times U$. More explicitly, we have the coherence equations
    \begin{equation}\label{eq:cocycle-condition}
    \beta_g\circ(g\cdot\beta_h)=\beta_{gh}\circ\alpha_{g,h,x},
    \qquad \beta_e\circ\iota_x=\id_x.
    \end{equation}
    A morphism $(x,\beta)\to(y,\beta')$ is an arrow $u:x\to y$ such that
    \[
    u\circ\beta_g=\beta'_g\circ(g\cdot u)
    \]
    as an identity over $G\times U$.
\end{definition}

Given descent data on an open cover $\{U_i\to U\}_i$, first glue the underlying objects in $\cX$. After transport to the glued object, the intertwining condition identifies the $\beta_i$ on $G\times(U_i\cap U_j)$. They glue to an isomorphism over $G\times U$. The cocycle and unit equations hold locally on $G^2\times U$ and $U$, respectively, and hence globally. The same sheaf condition glues intertwining arrows. For a prestack these formulas define the homotopy fixed prestack. Compare \cite[Definition~2.3 and Proposition~2.5]{romagny2005group}.

\begin{proposition}[The homotopy fixed locus biadjunction]\label{prop:hfix-adjunction}
    The homotopy fixed locus defines a pseudofunctor $\HomotopyFixedLocus:\St_G\to\St$, and there is a biadjunction
    \[
    \Trivial\dashv\HomotopyFixedLocus.
    \]
    Its unit and counit are the pseudonatural transformations
    \[
    \eta:\id_{\St}\Longrightarrow\HomotopyFixedLocus\circ\Trivial, \qquad \varepsilon:\Trivial\circ\HomotopyFixedLocus\Longrightarrow\id_{\St_G}.
    \]
    The unit sends $y$ to $(y,\id)$, and the counit has underlying morphism $(x,\beta)\mapsto x$, with equivariance constraint $\beta_g^{-1}:x\to g\cdot x$.
\end{proposition}
\begin{proof}
    For a $G$-morphism $(f,\sigma):\cX\to\cY$, set
    \[
    \HomotopyFixedLocus(f)(x,\beta)=
    \bigl(f(x),\{f(\beta_g)\circ\sigma_{g,x}^{-1}\}_{g\in G}\bigr).
    \]
    On arrows it is $u\mapsto f(u)$, and on equivariant $2$-morphisms it is the underlying natural transformation evaluated on the homotopy fixed objects.
    Write \(\delta_g=f(\beta_g)\circ\sigma_{g,x}^{-1}:g\cdot f(x)\to f(x)\).
    The constraints in \eqref{eq:cocycle-condition} have orientations \(\alpha_{g,h,x}:g\cdot(h\cdot x)\to(gh)\cdot x\) and \(\iota_x:x\to e\cdot x\). Naturality of \(\sigma_g\) at \(\beta_h:h\cdot x\to x\) gives
    \[
       \sigma_{g,x}^{-1}\circ (g\cdot f(\beta_h))
       = f(g\cdot\beta_h)\circ\sigma_{g,h\cdot x}^{-1}.
    \]
    Substitution, the cocycle for \(\beta\), and then the multiplication constraint of \((f,\sigma)\) give
    \[
    \begin{aligned}
       \delta_g\circ(g\cdot \delta_h) &=f(\beta_g)\circ \sigma_{g,x}^{-1}
            (g\cdot f(\beta_h))(g\cdot\sigma_{h,x})^{-1} = f\bigl(\beta_g(g\cdot\beta_h)\bigr) \circ \sigma_{g,h\cdot x}^{-1}\circ(g\cdot\sigma_{h,x})^{-1}\\
            &= f(\beta_{gh})\circ f(\alpha_{g,h,x})\circ\sigma_{g,h\cdot x}^{-1}(g\cdot\sigma_{h,x})^{-1} = f(\beta_{gh})\circ\sigma_{gh,x}^{-1}\circ\alpha'_{g,h,f(x)} \\
            &= \delta_{gh}\circ\alpha'_{g,h,f(x)}.
    \end{aligned}
    \]
    For the unit, \(\sigma_{e,x}\circ f(\iota_x)=\iota'_{f(x)}\) gives
    \[
       \delta_e\circ\iota'_{f(x)} = f(\beta_e)\circ \sigma_{e,x}^{-1} \circ \iota'_{f(x)} = f(\beta_e\circ\iota_x)=\id_{f(x)}.
    \]
    These are identities of smooth families over \(G\times G\times U\) and \(U\), respectively, so \(g,h\) also denote evaluation on arbitrary smooth sections after base change.

    For an arrow \(u:(x,\beta)\to(x',\beta')\), let \(\delta'_g\) be the transported structure on \(f(x')\). Then
    \[
       f(u)\circ\delta_g=f(\beta'_g)\circ f(g\cdot u)\circ\sigma_{g,x}^{-1} =\delta'_g(g\cdot f(u)).
    \]
    For an equivariant \(2\)-morphism \(\lambda:(f,\sigma)\Rightarrow(f',\sigma')\), naturality at \(\beta_g\) and equivariance at \(x\) give
    \[
       \lambda_x\circ\delta_g
       =f'(\beta_g)\circ\lambda_{g\cdot x}\circ\sigma_{g,x}^{-1}
       =f'(\beta_g)\circ(\sigma'_{g,x})^{-1}\circ(g\cdot\lambda_x).
    \]
    Finally, the constraint of a composite \((k,\tau)\circ(f,\sigma)\) at
    \((g,x)\) is \(\tau_{g,f(x)}\circ k(\sigma_{g,x})\). Transport along this
    composite gives \(k(\delta_g)\circ\tau_{g,f(x)}^{-1}\), precisely the result of
    the two successive transports. The identity morphism sends \(\beta_g\)
    to itself.
    The coherence constraints for composition and identity are those of the underlying morphisms, so this defines a pseudofunctor.

    The equivalence in \Cref{def:homotopy-fixed} sends a morphism $h:\cY\to\cX^{hG}$ to its underlying morphism equipped with the inverse fixing isomorphisms, i.e.,

    \[
    \begin{gathered}
    \Phi_{\cY,\cX}(h)=(\overline h,\sigma^h),\qquad
    h_U(y)=(x_y,\beta^y),\\
    \overline h_U(y)=x_y,\qquad
    \sigma^h_{g,y}=(\beta^y_g)^{-1}:x_y\longrightarrow g\cdot x_y.
    \end{gathered}
    \]
    The source action is trivial, so \(\overline h(g\cdot y)=x_y\).
    On \(v:y\to y'\), the map \(\overline h_U\) takes the underlying
    arrow of \(h_U(v)\). Its compatibility with \(\sigma^h\) is the
    intertwining equation for that arrow. The pullback constraints are
    those of \(h\) with their fixing data forgotten.
    Its inverse equips the values of an equivariant morphism from $\Trivial(\cY)$ with the inverse equivariance constraints. In detail, for \((k,\tau):\Trivial(\cY)\to\cX\), put
    \[
    \Psi_{\cY,\cX}(k,\tau)_U(y)
       =\bigl(k_U(y),\{\,\tau_{g,y}^{-1}\,\}_{g\in G(U)}\bigr),
    \qquad
    \Psi_{\cY,\cX}(k,\tau)_U(v)=k_U(v).
    \]
    The equivariance axioms reduce to
    \[
       \tau_{gh,y}=\alpha_{g,h,k(y)}\circ(g\cdot\tau_{h,y})\circ\tau_{g,y},
       \qquad \tau_{e,y}=\iota_{k(y)},
    \]
    and inverting these identities gives exactly the cocycle and unit for \(\tau^{-1}\). Naturality of \(\tau\) gives the arrow condition for \(k_U(v)\). 
    Equivariant $2$-morphisms correspond precisely to morphisms respecting the fixing data.
    Indeed, both functors take a \(2\)-cell to its underlying family of arrows. For \(\lambda:k\Rightarrow k'\) the required equivalence is
    \[
    (g\cdot\lambda_y)\circ\tau_{g,y}=\tau'_{g,y}\circ\lambda_y
       \quad\Longleftrightarrow\quad
    \lambda_y\circ\tau_{g,y}^{-1}=(\tau'_{g,y})^{-1}\circ(g\cdot\lambda_y).
    \]
    Applying \(\Phi\) after \(\Psi\) recovers \(k,\tau\) and applying \(\Psi\) after \(\Phi\) recovers \(h\), including its action on arrows and pullback constraints. In the displayed model of homotopy fixed locus the comparison isomorphisms have identity components.

    For \(v:\cY'\to\cY\) and an equivariant
    \((f,\sigma):\cX\to\cX'\), the naturality comparison is
    \[
      \Phi_{\cY',\cX'}\bigl(\HomotopyFixedLocus(f)\circ h\circ v\bigr)
       \ \cong\ (f,\sigma)\circ\Phi_{\cY,\cX}(h)\circ\Trivial(v).
    \]
    At \(y'\), put \(h(v(y'))=(x,\beta)\). Both sides have value \(f(x)\) with equivariance \(\sigma_{g,x}\circ f(\beta_g^{-1})\), and send an arrow to the same underlying composite. Thus identity components give this comparison in the explicit model. 
    The composite-transport formula above and associativity give compatibility with successive pre- and postcompositions. The calculation for \(\lambda\) gives compatibility with \(2\)-cells.

    Taking identity morphisms in the two Hom-groupoids gives the unit and counit
    \[
       \eta_{\cY}=\Psi_{\cY,\Trivial(\cY)}(\id_{\Trivial(\cY)}),\qquad
       \varepsilon_{\cX}=\Phi_{\HomotopyFixedLocus(\cX),\cX}(\id_{\HomotopyFixedLocus(\cX)}).
    \]
    The triangle modifications are
    \[\begin{tikzcd}
    	\Trivial & {\Trivial\HomotopyFixedLocus\Trivial} & \HomotopyFixedLocus & {\HomotopyFixedLocus\Trivial\HomotopyFixedLocus} \\
    	& \Trivial && \HomotopyFixedLocus
    	\arrow["{\Trivial\eta}", Rightarrow, nfold, from=1-1, to=1-2]
    	\arrow[""{name=0, anchor=center, inner sep=0}, "{\id_{\Trivial}}"', Rightarrow, nfold, from=1-1, to=2-2]
    	\arrow["{\varepsilon\Trivial}", Rightarrow, nfold, from=1-2, to=2-2]
    	\arrow["{\eta\HomotopyFixedLocus}", Rightarrow, nfold, from=1-3, to=1-4]
    	\arrow[""{name=1, anchor=center, inner sep=0}, "{\id_{\HomotopyFixedLocus}}"', Rightarrow, nfold, from=1-3, to=2-4]
    	\arrow["{\HomotopyFixedLocus\varepsilon}", Rightarrow, nfold, from=1-4, to=2-4]
    	\arrow["\Gamma"{pos=0.3}, between={0.2}{0.7}, Rightarrow, scaling nfold=3, from=1-2, to=0]
    	\arrow["\Delta"{pos=0.3}, between={0.2}{0.7}, Rightarrow, scaling nfold=3, from=1-4, to=1]
    \end{tikzcd}\]
    where, using this notation, we denote whiskering, so they have components
    \[
       \Gamma_{\cY}:\varepsilon_{\Trivial(\cY)}\circ \Trivial(\eta_{\cY})\Longrightarrow\id_{\Trivial(\cY)},\qquad \Delta_{\cX}:\HomotopyFixedLocus(\varepsilon_{\cX})\circ\eta_{\HomotopyFixedLocus(\cX)}\Longrightarrow\id_{\HomotopyFixedLocus(\cX)}.
    \]
    For $\Gamma$, the composite sends \(y\) to \((y,\id)\) and then to \(y\), sends an arrow \(v\) to \(v\), and has identity equivariance at every \(g\). Thus \((\Gamma_{\cY})_y=\id_y\) is equivariant. For $\Delta$, \(\eta_{H\cX}\) sends \((x,\beta)\) to \(((x,\beta),\id)\). Transporting this outer identity by the counit gives the fixing arrow
    \[
       \varepsilon_{\cX}(\id_{(x,\beta)})\,(\beta_g^{-1})^{-1}=\beta_g.
    \]
    Thus \((\Delta_{\cX})_{(x,\beta)}=\id_x:(x,\beta)\to(x,\beta)\), and the composite preserves each intertwining arrow. These identity components commute with every arrow and every change of parameter.
    Both modifications are invertible, and their two coherence pastings are identities in the explicit model.
    In conclusion, in the specific model of the homotopy fixed locus with pairs $(x,\beta)$, we obtain a strict biadjunction; transport along the canonical comparison isomorphisms in an equivalent model gives the same biadjunction, not necessarily strict.
\end{proof}

\begin{corollary}\label{cor:hfix-limits}
    The homotopy fixed locus preserves all existing small bilimits in $\St_G$ (with values in $\St$). In particular, for a diagram $\cX\to\mathcal Z\leftarrow\cY$ of $G$-stacks,
    \[
    \HomotopyFixedLocus(\cX\times_{\mathcal Z}\cY) \simeq\HomotopyFixedLocus(\cX)\times_{\HomotopyFixedLocus(\mathcal Z)}\HomotopyFixedLocus(\cY),
    \]
    where the left-hand fiber product carries the induced coherent action. Also, $\pt^{hG}\simeq\pt$ for the trivial action.
\end{corollary}
\begin{proof}
    As a right biadjoint, $\HomotopyFixedLocus$ preserves bilimits.
\end{proof}

\begin{remark}[The trivial-action case]\label{rem:hfix-trivial}
    For every stack $\cY$ there is a canonical equivalence
    \[
    \HomotopyFixedLocus\Trivial(\cY)\simeq\underline{\Hom}(\BG,\cY).
    \]
    Indeed, descent along $U\to\BG\times U$ identifies a morphism $\BG\times U\to\cY$ with an object $y\in\cY(U)$ and a smooth homomorphism of group sheaves $G_U\to\Aut_{\cY}(y)$, with arrows intertwining the homomorphisms. Under this identification the unit of \Cref{prop:hfix-adjunction} sends an object to the constant family, i.e. the trivial homomorphism, and evaluation along $\pt\to\BG$ is a retraction of that unit. If $\Trivial$ also had $\FixedLocus$ as right biadjoint, uniqueness of a right biadjoint would give a pseudonatural equivalence $\FixedLocus\simeq\HomotopyFixedLocus$. On the trivial $\mathbb Z/2$-action on $\mathcal B(\mathbb Z/2)$ one side has one isomorphism class over $\pt$, the other has two, by \Cref{ex:trivial-fixed-action}.
\end{remark}

\begin{remark}
    Over an object $x:U\to\cX$, the fiber of $\cX^{hG}\to\cX$ is the sheaf of group-theoretic sections of $\alpha_x:\Stab^G_{\cX}(x)\to G_U$, i.e. the group-theoretic Weil restriction of $\Stab_{\cX}^G$ along the group stack $G\times\cX\to\cX$. This is the smooth counterpart of \cite[Lemma~4.1.2]{romagny2022algebraicity}.
\end{remark}

\begin{lemma}\label{lemma:action}
    A smooth action on a category fibred on groupoids $\cX$ extends to its stackification $a:\cX\to\cX^{\st}$, and $a$ is equivariant. The extension is canonical up to coherent equivalence, even when the action on the chosen model is strict.
\end{lemma}
\begin{proof}
    Stackification preserves finite $2$-limits and $G$ is already a sheaf, so $(G\times\cX)^{\st}\simeq G\times\cX^{\st}$. Extend the action morphism and its unit and associativity cells by the universal property of stackification (see for example \cite[Theorem 3.5.18]{alper2025stacks} for the stackification construction). The pentagon and unit equations are preserved, since they are equalities between the induced $2$-morphisms. The uniqueness is bicategorical uniqueness.
\end{proof}

\begin{proposition}
\label{prop:stackificationfixed}
    If $\cX$ is a prestack, then the canonical morphism
    \[
    (\cX^{hG})^{\st}\longrightarrow(\cX^{\st})^{hG}
    \]
    is an isomorphism of stacks.
\end{proposition}
\begin{proof}
    The map $a:\cX\to\cX^{\st}$ is fully faithful and every object of the target locally comes from $\cX$. Choose an open cover of $U$ on which the underlying object of $(x,\beta)$ lifts. On each member $V$ of the cover, full faithfulness on $G\times V$ lifts the entire isomorphism $\beta$. Faithfulness on $G^2\times V$ and $V$ gives its cocycle and unit equations. The same argument lifts morphisms. Thus $\cX^{hG}\to(\cX^{\st})^{hG}$ is fully faithful and locally essentially surjective, which proves the claim after stackification.
\end{proof}

\begin{corollary}[Stackification and the four functors]\label{cor:fixed-stackification}
    For a prestack $\cY$ and a $G$-prestack $\cX$, the canonical comparisons are equivalences for all four functors with stackification:
    \[
    \begin{aligned}
    \Trivial(\cY)^{\st}&\simeq\Trivial(\cY^{\st}),&
    \Forgetful(\cX)^{\st}&\simeq\Forgetful(\cX^{\st}),\\
    \FixedLocus(\cX)^{\st}&\simeq\FixedLocus(\cX^{\st}),&
    \HomotopyFixedLocus(\cX)^{\st}&\simeq\HomotopyFixedLocus(\cX^{\st}).
    \end{aligned}
    \]
\end{corollary}
\begin{proof}
    The stackification morphism $a:\cX\to\cX^{\st}$ is fully faithful for a prestack and equivariant by \Cref{lemma:action}. For any $x\in\cX(U)$ it induces an isomorphism of extended stabilizer sheaves over $G_U$: full faithfulness identifies every fixing isomorphism after pullback to any $V\to U$. Hence $x$ is fixed if and only if $a(x)$ is fixed. The restriction to the fixed subprestacks is fully faithful and locally essentially surjective, since every object of $\cX^{\st}$ locally lifts to $\cX$, and the fixed condition holds for each such lift. Its stackification is therefore an equivalence. The homotopy-fixed comparison is \Cref{prop:stackificationfixed}. The trivial and forgetful comparisons follow from the construction of the equivariant stackification, whose underlying stack is the ordinary stackification and which sends the trivial constraints to the canonical trivial constraints. All comparisons arise from the equivariant stackification morphism and are coherent with the restriction and forgetful maps.
\end{proof}

There is a canonical faithful forgetful morphism
\[
\epsilon_{\cX}:\cX^{hG}\longrightarrow\cX,\qquad(x,\beta)\longmapsto x.
\]
\begin{proposition}\label{prop:hfp-factorization}
    It factors canonically through $\cX^G$:
    \[
    \begin{tikzcd}
    {\cX^{hG}} & \cX \\
    & {\cX^G}.
    \arrow["\epsilon_{\cX}",from=1-1,to=1-2]
    \arrow["j_{\cX}"',from=1-1,to=2-2]
    \arrow["i_{\cX}"',hook,from=2-2,to=1-2]
    \end{tikzcd}
    \]
\end{proposition}
\begin{proof}
    For $u:V\to U$ and $g:V\to G$, pull $\beta$ back along $(g,u):V\to G\times U$. The resulting arrow $g\cdot u^*x\to u^*x$ gives a section of $\alpha_x$ as a homomorphism of group sheaves by the cocycle and unit equations. Hence $\alpha_x$ is an epimorphism. See also \cite[Remark~A.17]{aranha2025virtual}.
\end{proof}

\begin{proposition}[Pseudonatural comparison]\label{prop:fixed-pseudonatural}
    There are canonical pseudonatural transformations
    \[
    j:\HomotopyFixedLocus\Longrightarrow\FixedLocus,
    \qquad i:\FixedLocus\Longrightarrow\Forgetful,
    \qquad \epsilon:\HomotopyFixedLocus\Longrightarrow\Forgetful.
    \]
    To describe $\epsilon$ in terms of the counit of \Cref{prop:hfix-adjunction}, write $\phi:\Forgetful\Trivial\Rightarrow\id_{\St}$ for the canonical identification. Choose a quasi-inverse $\phi^{-1}$ with its unit and counit modifications. Then $\epsilon$ is the composite
    \[
           \epsilon:\HomotopyFixedLocus \xRightarrow{\,\phi^{-1}\HomotopyFixedLocus\,} \Forgetful \Trivial\HomotopyFixedLocus
             \xRightarrow{\,\Forgetful\varepsilon\,}\Forgetful.
    \]
    In the explicit model \(\Forgetful\Trivial\) forgets a trivial action, so the component at \(\cX\) sends \((x,\beta)\) to \(x\) and \(u\) to \(u\), and the equivariant counit also carries the arrow \(\beta_g^{-1}\), which is forgotten upon applying \(\Forgetful\).
    There is a canonical invertible modification
    \[\begin{tikzcd}
    	\HomotopyFixedLocus && \FixedLocus \\
    	& \Forgetful
    	\arrow["j", Rightarrow, nfold, from=1-1, to=1-3]
    	\arrow[""{name=0, anchor=center, inner sep=0}, "{\epsilon}"', Rightarrow, nfold, from=1-1, to=2-2]
    	\arrow["i", Rightarrow, nfold, from=1-3, to=2-2]
    	\arrow["\Lambda"{pos=0.75},between={0.5}{0.7}, Rightarrow, scaling nfold=3, from=1-3, to=0]
    \end{tikzcd}\]
    such that in the explicit models it is the identity $i\circ j=\epsilon$.
\end{proposition}
\begin{proof}
    The factorization on objects is \Cref{prop:hfp-factorization}. Given a $G$-morphism $(f,\sigma)$, both composites in the naturality square for $j$ send $(x,\beta)$ to $f(x)$, and both send an arrow $u$ to $f(u)$. The coherence for composition follows from the underlying composition constraints. Naturality with respect to equivariant $2$-morphisms is likewise inherited from the underlying $2$-morphisms. The comparison with the counit follows from the formula for $\varepsilon$ in \Cref{prop:hfix-adjunction}.
    Write \(k:\mathrm{P}\Rightarrow \mathrm{Q}\) for any of \(j,i,\epsilon\). For each
    \(a=(f,\sigma):\cX\to\cY\), the naturality square comes with an
    invertible \(2\)-cell in \(\St\) of type
    \[
        k_a:\mathrm{Q}(a)k_{\cX}\Longrightarrow k_{\cY}\mathrm{P}(a).
    \]
    In the three cases its source and target are
    \[
    \begin{aligned}
       j_a &: \FixedLocus(a) j_{\cX}\Longrightarrow j_{\cY}\HomotopyFixedLocus(a),\\
       i_a &: \Forgetful(a)i_{\cX}\Longrightarrow i_{\cY}\FixedLocus(a),\\
       \epsilon_a &: \Forgetful(a)\epsilon_{\cX}\Longrightarrow\epsilon_{\cY}\HomotopyFixedLocus(a).
    \end{aligned}
    \]
    For example, the naturality cell of the equivariant counit has
    underlying component \(\id_{f(x)}\). Its equivariance follows from
    \[
       \sigma_{g,x}\circ f(\beta_g^{-1})
       =\bigl(f(\beta_g)\circ\sigma_{g,x}^{-1}\bigr)^{-1}.
    \]
    Its image under \(\Forgetful\) gives \(\epsilon_a\) with the identification
    \(\phi\) above. The comparison cells for \(i,j\) also have identity
    components in the model of full substacks and pairs.

    For composable \(a,b\), the cell for \(ba\) is the pasting of
    \(\mathrm{Q}(b)k_a\) and \(k_b \mathrm{P}(a)\), using the compositors of \(\mathrm{P},\mathrm{Q}\) to identify its endpoints. The identity axiom uses their unitors.
    For a \(2\)-cell \(\lambda:a\Rightarrow a'\) the required equality is
    \[
       (k_{\cY}\mathrm{P}(\lambda))\circ k_a
       =k_{a'}\circ(\mathrm{Q}(\lambda)k_{\cX}).
    \]
    Here all three operations send \(\lambda\) to its underlying components at the indicated objects. This proves the equality.
    Pseudonaturality records precisely these invertible comparison cells and equations.

    Regarding the modification, at each \(\cX\) it gives a \(2\)-isomorphism \(\Lambda_{\cX}:i_{\cX}j_{\cX}\Rightarrow\epsilon_{\cX}\).
    If \((ij)_a\) is the pasting of \(i_a j_{\cX}\) and \(i_{\cY}j_a\), its compatibility at \(a\) reads
    \[
       (\Lambda_{\cY}\HomotopyFixedLocus(a))\circ(ij)_a=\epsilon_a\circ(\Forgetful(a)\Lambda_{\cX}).
    \]
    At \((x,\beta)\) both components of the triangle forget \(\beta\), so \(\Lambda_{\cX}(x,\beta)=\id_x\) in the explicit model, and both sides of the last display evaluate to \(\id_{f(x)}\) there; in particular, equality can hold with these choices. The displayed modification specifies the comparison that is transported when models are changed. This proves the pseudonaturality and the identity of the two composites in the chosen models.
\end{proof}

\begin{remark}\label{rem:comparison-fibers}
For $x:U\to\Forgetful(\cX)$ there is a canonical isomorphism of stacks over $U$
\[
U\times_{\Forgetful(\cX)}\HomotopyFixedLocus(\cX)
\simeq\underline{\operatorname{Spl}}(\alpha_x),
\]
where the right-hand side is the sheaf on $\Man/U$ of group-sheaf sections of $\alpha_x:\Stab^G_{\cX}(x)\to G_U$. Explicitly, over $u:V\to U$ it consists of homomorphisms $s:G_V\to\Stab^G_{\cX}(u^*x)$ with $\alpha_{u^*x}\circ s=\id_{G_V}$. An isomorphism to the specified underlying object transports any fixing structure to that object, and the remaining morphisms over its identity are identities. This proves the equivalence as sheaves. If $x$ is fixed, the same sheaf is the fiber of $j_{\cX}$ along $x:U\to\cX^G$, since $i_{\cX}$ is fully faithful. Thus $j_{\cX}$ and $\epsilon_{\cX}$ are faithful and representable by sheaves of sets, but need not be full, epimorphic, or representable by manifolds. The inclusion $i_{\cX}$ is fully faithful. In particular, the condition that a fixed object be in the local essential image of $j_{\cX}$ is the existence of a splitting locally on its test manifold; it is stronger than the epimorphism condition defining the fixed locus. The non-split extension in \Cref{ex:trivial-fixed-action} illustrates this distinction.
\end{remark}

For a differentiable $G$-stack put $\cN=[\cM/G]$, and let $r:\cN\to\BG$ classify its principal $G$-bundle $q:\cM\to\cN$. There are natural equivalences
\begin{equation}\label{eq:hfp-section-stack}
\cM^{hG}\simeq\underline{\Hom}_G(\pt,\cM)
\simeq\underline{\Hom}_{\BG}(\BG,\cN).
\end{equation}
On $T$ the last groupoid consists of $s:\BG\times T\to\cN$ with a specified frame $\theta:r s\Rightarrow\operatorname{pr}_{\BG}$. An arrow $\eta:(s,\theta)\to(s',\theta')$ satisfies $\theta'\circ(r\eta)=\theta$. For the corresponding algebraic description see \cite[Remark~A.16]{aranha2025virtual}.

\begin{proposition}\label{prop:quotient-adjunction}
    For stacks, the quotient construction is left biadjoint to the trivial-action functor, $[-/G]\dashv\Trivial$, i.e.
    \[
    \Mor_{\St}([\cX/G],\cY)
    \simeq\Mor_{\St_G}(\cX,\Trivial(\cY)).
    \]
    In particular, the trivial action satisfies $[\Trivial(\cY)/G]\simeq\BG\times\cY$ over $\BG$. For a differentiable $G$-stack $\cM$, the equivalence \eqref{eq:hfp-section-stack} identifies
    \[
    \cM^{hG}\simeq\underline{\Hom}_{\BG}(\BG,[\cM/G])
    \]
    pseudonaturally in $\cM$. Under this equivalence, $\epsilon_{\cM}$ is the restriction of a section along $\pt\to\BG$, followed by the canonical identification $[\cM/G]\times_{\BG}\pt\simeq\cM$.
\end{proposition}
\begin{proof}
    A morphism from the quotient to $\cY$ is, by descent along the quotient map, a morphism from $\cX$ equipped with compatible isomorphisms identifying its translates. These are exactly the equivariance data for a morphism into the trivial action, with the same cocycle, unit and arrow conditions. The correspondence is compatible with composition and $2$-morphisms, giving the biadjunction. For the trivial action the quotient description on a test manifold consists of a principal $G$-bundle and an independent object of $\cY$, which gives $\BG\times\cY$ with its projection to $\BG$.

    For \eqref{eq:hfp-section-stack}, use the atlas $T\to\BG\times T$ determined by the trivial torsor, with groupoid $G\times T\rightrightarrows T$. Together, a section restricts to $y\in\cN(T)$ and the frame $r(y)\simeq T\times G$ are an object $x$ of $\cN\times_{\BG}\pt\simeq\cM$ over $T$. Over $G\times T$ the descent arrow is an automorphism of the pullback of $y$. Compatibility with $\theta$ says that on the framed torsor its value at $g$ is left multiplication by $g$. In the transformation prestack a degree-$g$ arrow of $q(x)$ is precisely an isomorphism $\beta_g:g\cdot x\to x$ \emph{over $G\times T$}. This description also holds after stackification, since the degree and the arrow are sections of sheaves. Over $G^2\times T$ the composition law of the transformation prestack is exactly \eqref{eq:cocycle-condition}; the identity gives its unit equation.
    
    Conversely a smooth $\beta$ satisfying these two equations gives effective descent for $q(x)$ and for its frame, hence a section $(s,\theta)$. An arrow of framed descent data has degree $e$ on $T$ and intertwines $\beta$ on $G\times T$, precisely the arrow condition in $\cM^{hG}(T)$. Pullback along $T'\to T$ pulls back all these data along $G^n\times T'\to G^n\times T$ for $n=0,1,2$. This gives the pseudonatural inverse equivalences and identifies evaluation with forgetting the descent arrow.
\end{proof}

For the chosen groupoid presentation, a homotopy fixed structure on $x\in X_0$ is exactly a smooth group-theoretic section $\sigma(g)=(g,\beta_g)$ of $\pi_x:\widehat G_x\to G$. Thus
\[
\begin{aligned}
\text{fixed locus at }x&:\quad\pi_x\text{ is surjective},\\
\text{homotopy fixed locus over }x&:\quad\pi_x\text{ has a chosen smooth homomorphic section}.
\end{aligned}
\]

\subsection{Quotient stacks}

Let $K$ be a Lie group and $\qman$ be a $K$-manifold. A $G$-action on $[\qman/K]$ need not lift to the chosen atlas, we need to “extend'' it. We study the actions on this presentation induced by a Lie group extension
\begin{equation}\label{eq:group-extension}
\begin{tikzcd}
1 & K & {\widetilde G} & G & 1,
\arrow[from=1-1,to=1-2]
\arrow[from=1-2,to=1-3]
\arrow["\pi",from=1-3,to=1-4]
\arrow[from=1-4,to=1-5]
\end{tikzcd}
\end{equation}
where $\pi$ is a surjective submersion, together with a smooth left $\widetilde G$-action on $\qman$, denoted $\widetilde g\ast m$. It induces a $G$-action on $[\qman/K]$. One precise construction is the identification
\begin{equation}\label{eq:extension-fibre}
[\qman/K]\simeq[\qman/\widetilde G]\times_{\BG}\pt,
\qquad [[\qman/K]/G]\simeq[\qman/\widetilde G],
\end{equation}
where the morphism to $\BG$ is induced by $\pi$. As in the last step of \Cref{lem:relative-slice}, the first fiber is $[(\qman\times G)/\widetilde G]$ with action $(m,a)\mapsto(\widetilde g\ast m,\pi(\widetilde g)a)$; the reduction to $a=e$ gives $[\qman/K]$. The principal bundle action on this fiber defines the residual action and the second equivalence.

Write
\[
S_\pi:=\operatorname{spl}(\pi)
=\Set{\sigma\in\hom(G,\widetilde G)|\pi\circ\sigma=\id_G},
\quad K_\sigma:=C_K(\sigma)=K\cap C_{\widetilde G}(\sigma).
\]

\begin{lemma}\label{lem:splitting-space}
    If $G$ is compact, then $S_\pi$ is a countable disjoint union of open $K$-orbits, each diffeomorphic to $K/K_\sigma$. It represents the sheaf of smooth families of splittings. A smooth family of splittings is locally $\sigma_t=k(t)\sigma_{t_0}k(t)^{-1}$ for a smooth $K$-valued map $k$.
\end{lemma}
\begin{proof}
    Use \Cref{thm:homspace} with target $\widetilde G$. If $\sigma'=h\sigma h^{-1}$ and both homomorphisms split $\pi$, then for every $g\in G$
    \[
    g=\pi(h)g\pi(h)^{-1},\qquad\pi(h)\in Z(G).
    \]
    Consequently
    \[
    k=h\sigma(\pi(h))^{-1}\in K,\qquad
    \sigma'=k\sigma k^{-1},
    \]
    because $\sigma(\pi(h))$ centralizes $\sigma(G)$. Thus the intersection of any $\widetilde G$-conjugacy orbit with $S_\pi$ is either empty or one $K$-orbit. It is open, and there are countably many such orbits. Their stabilizers are the closed groups $K_\sigma$.
    
    For a smooth family, \Cref{thm:homspace} provides smooth conjugators $h(t)\in\widetilde G$. The same formula $k(t)=h(t)\sigma_{t_0}(\pi(h(t)))^{-1}$ gives smooth conjugators in $K$. It also gives the continuous local inverses for the homogeneous orbit topology. Conversely, local lifts to $K$ make evaluation smooth. This proves the asserted manifold and sheaf statements.
\end{proof}

\begin{theorem}
\label{thm:homotopy-fixed-quotient}
    In \eqref{eq:group-extension}, assume $G$ is compact and let $\qman$ be a $\widetilde G$-manifold. Give $[\qman/K]$ the action induced by \eqref{eq:extension-fibre}. For $\sigma\in S_\pi$ put
    \[
    \qman^\sigma=\Set{m\in \qman | \sigma(g)\ast m=m \ \forall\,g\in G},
    \]
    and define the $K$-set
    \[
    Z=\Set{(m,\sigma)\in \qman\times S_\pi | m\in \qman^\sigma},\qquad
    k\cdot(m,\sigma)=(k\ast m,k\sigma k^{-1}).
    \]
    Then $S_\pi$ and $Z$ are manifolds, and there are isomorphisms of stacks
    \begin{equation}\label{eq:hfp-quotient-decomposition}
    [\qman/K]^{hG}\simeq[Z/K]
    \simeq\coprod_{[\sigma]\in S_\pi/K}[\qman^\sigma/K_\sigma].
    \end{equation}
    In particular $[\qman/K]^{hG}$ is a differentiable stack. The map $\epsilon_{[M/K]}$ is induced by $(m,\sigma)\mapsto m$, and
    \begin{equation}\label{eq:epsilon-atlas}
    \qman\times_{[\qman/K]}[\qman/K]^{hG}\simeq Z.
    \end{equation}
    Thus $\epsilon_{[M/K]}$ is quasi-representable. If $S_\pi=\varnothing$, the homotopy fixed locus is empty.
\end{theorem}
\begin{proof}
    \emph{1. The equivalence on smooth families.}
    Using \eqref{eq:hfp-section-stack} and \eqref{eq:extension-fibre}, take a framed section over $\BG\times T$. Descent along $T\to\BG\times T$ gives a right $\widetilde G$-torsor $Q\to T$, a map $f:Q\to \qman$ with $f(qh)=h^{-1}\ast f(q)$, and a commuting left $G$-action on $Q$ preserving $f$ and covering $\id_T$. To see the smoothness, the descent arrow is a bundle isomorphism over $G\times T$. Under the canonical pullback identifications its total-space map is $(g,q)\mapsto(g,g\cdot q)$, which is smooth in both variables. Descent over $G^2\times T$ and along the unit gives the action equations. The specified $2$-cell is a trivialization $Q\times_{\widetilde G}G\simeq T\times G$ intertwining left multiplication by $G$. Equivalently it is a smooth function $b:Q\to G$ satisfying
    \[
                 b(qh)=b(q)\pi(h),\qquad b(g\cdot q)=g\,b(q).
    \]
    
    The normalized frames form the right $K$-torsor $P=b^{-1}(e)\to T$. Indeed any local section $q$ of $Q$ can be normalized by multiplying it by a local smooth lift of $b(q)^{-1}$ through $\pi$. For $p\in P$ define
    \[
        m(p)=f(p),\qquad g\cdot p=p\,\sigma_p(g).
    \]
    The division map of the $\widetilde G$-torsor gives joint smoothness in $(p,g)$. Commutation of the actions makes $\sigma_p$ a homomorphism, the frame gives $\pi\sigma_p=\id_G$, and $f$-invariance gives $\sigma_p(g)\ast m(p)=m(p)$. Changing frame gives, for $c\in K$,
    \[
        m(pc)=c^{-1}\ast m(p),\qquad \sigma_{pc}=c^{-1}\sigma_p c.
    \]
    By \Cref{lem:splitting-space} and the embedded description of $Z$ below, these are precisely the data of a smooth $K$-equivariant map $P\to Z$ in the right-torsor convention.
    
    Conversely, such a $P\to Z$ gives $Q=P\times_K\widetilde G$, with
    \[
      g\cdot[p,h]=[p,\sigma_p(g)h],\qquad
      f([p,h])=h^{-1}\ast m(p),\qquad b([p,h])=\pi(h).
    \]
    The change-of-frame equations make these formulas well-defined and smooth. They supply the descent arrow, its cocycle and the frame. A morphism of sections restricts to a map of their normalized torsors preserving $(m,\sigma)$, and any such map extends by $[p,h]\mapsto[\psi(p),h]$. Multiplication $[p,h]\mapsto ph$ and the inclusion $p\mapsto[p,e]$ identify the two composites with the identity. These bundle isomorphisms commute with pullback in $T$, giving the first equivalence in \eqref{eq:hfp-quotient-decomposition}.
    
    \emph{2. Smoothness and the orbit decomposition.}
    By \Cref{lem:splitting-space}, $S_\pi$ is a manifold with smooth evaluation. We have a smooth $G$-action on $\qman\times S_\pi$
    \[
    g\cdot(m,\sigma)=(\sigma(g)\ast m,\sigma)
    \]
    whose fixed set is $Z$. A $G$-invariant Riemannian metric on the product gives, at any fixed point, exponential coordinates identifying that fixed set with the fixed tangent subspace. Thus $Z$ is embedded and is closed by the Hausdorff property. The same argument applies to $\qman^\sigma$.
    
    If $c\in K_\sigma$ and $m\in \qman^\sigma$, then
    \[
    \sigma(g)\ast(c\ast m)=c\ast(\sigma(g)\ast m)=c\ast m,
    \]
    so $K_\sigma$ acts on $\qman^\sigma$. For the open orbit $\mathcal O_\sigma\subseteq S_\pi$, the corresponding open part of $Z$ is $K$-equivariantly diffeomorphic to
    \begin{equation}\label{eq:associated-fixed-piece}
    K\times_{K_\sigma}\qman^\sigma\longrightarrow Z_{[\sigma]},
    \qquad[k,m]\longmapsto(k\ast m,k\sigma k^{-1}).
    \end{equation}
    Here $(k,m)\sim(kc,c^{-1}\ast m)$makes this map well defined, and the map is surjective. Equality of two images implies $k'^{-1}k\in K_\sigma$ and the corresponding relation on the points, so it is injective. A local smooth section of $K\to K/K_\sigma$ gives its smooth inverse. The action of $K_\sigma$ on $K\times \qman^\sigma$ is free and proper, because its right action on $K$ is proper for a closed subgroup; hence the associated bundle is a manifold.
    
    The inclusion $K_\sigma\ltimes \qman^\sigma\to K\ltimes Z_{[\sigma]}$ is fully faithful: an arrow between $(m,\sigma)$ and $(m',\sigma)$ belongs precisely to $K_\sigma$. It is essentially surjective by the surjective submersion $K\times \qman^\sigma\to Z_{[\sigma]}$. It is therefore a Morita equivalence and gives $[Z_{[\sigma]}/K]\simeq[\qman^\sigma/K_\sigma]$. The countable disjoint union gives the second equivalence and a second-countable atlas.
    
    \emph{3. The forgetful map.}
    Under the first equivalence, forgetting the section data is induced by the $K$-equivariant projection $Z\to \qman$. The standard pullback square for this map of quotient stacks gives \eqref{eq:epsilon-atlas}, including the chosen trivialization of the underlying principal $K$-bundle. This proves quasi-representability. More generally, a map $T\to[\qman/K]$ is a principal $K$-bundle $P\to T$ with equivariant map $f:P\to \qman$. Its base change of $\epsilon$ is the quotient sheaf
    \[
    (P\times_\qman Z)/K,
    \]
    with the fiber product first formed in sheaves. For instance take $K=1$, $G=\mathbb Z/2$ acting by sign on $\qman=\mathbb R$, and test $0\hookrightarrow\mathbb R$ along $(a,b)\mapsto ab$. Then the resulting sheaf is not representable by a manifold.
\end{proof}

\begin{corollary}\label{cor:semidirect}
    Suppose that $G$ is compact and acts smoothly on $\qman$ and on $K$, with the action $\theta$ on $K$ by group automorphisms. Assume the compatibility condition
    \[
    g\cdot(k\ast m)=\theta_g(k)\ast(g\cdot m).
    \]
    Then $\widetilde G=K\rtimes_\theta G$ acts by $(k,g)\ast m=k\ast(g\cdot m)$ and
    \[
    [\qman/K]^{hG}\simeq
    \coprod_{[\varphi]\in Z^1(G,K)/K}[\qman^\varphi/K_\varphi],
    \]
    where $Z^1(G,K)=\Set{\varphi\in C^\infty(G,K)| \varphi_{gh}=\varphi_g\theta_g(\varphi_h)\ \forall\,g,h\in G}$ is the set of non-abelian 1-cocycles, with the $K$-action defined as a $\theta$-twisted conjugation
    \[
    (k\cdot\varphi)_g=k\varphi_g\theta_g(k)^{-1},
    \]
    and $\qman^\varphi=\Set{m\in \qman|\varphi_g\ast(g\cdot m)=m \ \forall\, g\in G}$, with an action given by the isotropy $K_\varphi=\Set{k\in K|k\varphi_g\theta_g(k)^{-1}=\varphi_g \ \forall\, g\in G}$.
\end{corollary}
\begin{proof}
    The multiplication is $(k,g)(l,h)=(k\theta_g(l),gh)$. The correspondence
    \[
    \varphi\longleftrightarrow\widetilde\varphi,
    \qquad\widetilde\varphi(g)=(\varphi_g,g)
    \]
    identifies cocycles with $\operatorname{spl}(\pi)$, because
    \[
    \widetilde\varphi(g)\widetilde\varphi(h)
    =(\varphi_g\theta_g(\varphi_h),gh)
    =\widetilde\varphi(gh).
    \]
    Conjugation by $(k,e)$ gives
    \[
    (k,e)(\varphi_g,g)(k,e)^{-1}
    =(k\varphi_g\theta_g(k)^{-1},g).
    \]
    Hence the gauge stabilizer is exactly
    \[
    K_\varphi=K\cap C_{\widetilde G}(\widetilde\varphi),
    \]
    and $\qman^\varphi=\qman^{\widetilde\varphi}$. By \Cref{lem:splitting-space} and \Cref{thm:homotopy-fixed-quotient}, $K_\varphi$ preserves $\qman^\varphi$, and the piece of the parameter manifold over the gauge orbit is $K\times_{K_\varphi}\qman^\varphi$.
\end{proof}

\begin{remark}\label{rem:direct}
    If $\theta_g=\id_K$, then $Z^1(G,K)=\hom(G,K)$, $K_\varphi=C_K(\varphi)$, and the formula becomes
    \[
    [\qman/K]^{hG}\simeq
    \coprod_{[\varphi]\in\hom(G,K)/K}[\qman^\varphi/C_K(\varphi)].
    \]
    For the associated bundles with compact target, see \cite[Lemma~17 and Proposition~18]{Mrcun2018}. The proof above also identifies the framed section stack on smooth test families.
    
    If the action of $G=\mathbb Z/k$ on the stack is trivial, a homomorphism is determined by the image $a$ of the specified generator. Consequently
    \[
     \underline{\Hom}(\mathcal B(\mathbb Z/k),[\qman/K])
     \simeq\coprod_{\substack{[a]\ \mathrm{in}\ K\\a^k=e}} [\qman^a/C_K(a)].
    \]
    For an orbifold this is the part of the inertia stack whose element has order dividing $k$, including the untwisted sector. Restricting to injections $\mathbb Z/k\hookrightarrow K_m$ selects the exact-order-$k$ sectors. This qualification is needed for the injectivity convention in \cite[\S~5.1]{jiang2026equivariant}. The fixed-subspace and centralizer description is the classical one of \cite[\S~3.1, Example~3.1.3]{chenruan}.

    \Cref{thm:homotopy-fixed-quotient} also applies to nonsplit extensions. In that case the atlas $\qman\to[\qman/K]$ may fail to be a $G$-atlas.
\end{remark}

\subsection{Separated stacks}

Romagny proves algebraicity and smoothness results for algebraic fixed point stacks under the hypotheses of \cite[Theorems~4.1.4 and~4.3.6]{romagny2022algebraicity}. On the smooth site we use compactness and the relative slice of \Cref{lem:relative-slice} to prove differentiability of the homotopy fixed locus and quasi-representability of evaluation.

\begin{theorem}[Differentiability of the homotopy fixed locus]
\label{thm:hfp-differentiability}
    Let $G$ be a compact Lie group acting smoothly on a separated differentiable stack $\cM$. Then $\cM^{hG}$ is a separated differentiable stack and $\epsilon_{\cM}:\cM^{hG}\to\cM$ is quasi-representable. At a point $(x,\sigma)$, write $H_x=\Aut_{\cM}(x)$, the isotropy group is
    \[
        \Aut_{\cM^{hG}}(x,\sigma)
           \cong H_x\cap C_{\widehat G_x}(\sigma).
    \]
    
    For each homotopy fixed object $(x,\sigma_0)$ over a point there is an open $G$-invariant neighborhood $\cW_x\subseteq\cM$ of $x$ and a $G$-equivariant equivalence $\cW_x\simeq[V_x/H_x]$ with action induced by the compact extension
    \[\begin{tikzcd}
    	1 & {H_x} & {\widehat G_x} & G & 1.
    	\arrow[from=1-1, to=1-2]
    	\arrow[from=1-2, to=1-3]
    	\arrow["{\pi_x}", from=1-3, to=1-4]
    	\arrow[from=1-4, to=1-5]
    \end{tikzcd}\]
    Here $V_x$ is an invariant open ball in the normal representation at $x$, and the corresponding open substack is
    \begin{equation}\label{eq:hfp-local-model}
    \cW_x^{hG}\simeq
    \coprod_{[\sigma]\in\operatorname{spl}(\pi_x)/H_x}
    [\,V_x^\sigma/C_{H_x}(\sigma)\,]
    \subseteq\cM^{hG}.
    \end{equation}
    The sector containing $(x,\sigma_0)$ is the one indexed by $[\sigma_0]$.
\end{theorem}
\begin{proof}
    \emph{1. The quotient over $\BG$.}
    Put $\cN=[\cM/G]$ and let $r:\cN\to\BG$ classify the principal bundle $q:\cM\to\cN$. By \Cref{prop:equivariant-atlas,prop:proper-quotient}, $\cN$ is separated and differentiable, and by \eqref{eq:hfp-section-stack}, our target is the section stack of $r$.

    \emph{2. The relative slice.}
    For $y=q(x)$, its isotropy is $\widehat G_x$ and the homomorphism induced by $r$, using the lift $x$ as framing, is $\pi_x$. The given homotopy fixed structure splits $\pi_x$, so it is surjective and its kernel is $H_x$. By \Cref{lem:relative-slice} we have an open chart
    \[
    [V_x/\widehat G_x]\simeq\cU_x\subseteq\cN
    \]
    and a specified $2$-isomorphism between $r|_{\cU_x}$ and the map induced by $\pi_x$. The essential step is the equivariant trivialization of the principal $G$-bundle over $V_x$: averaging a connection over the compact group $\widehat G_x$ and transporting radially gives
    \[
    P_x\cong V_x\times G,\qquad
    \widehat g\cdot(v,a)=(\widehat g\cdot v,\pi_x(\widehat g)a).
    \]
    Consequently, the open pullback $\cW_x=q^{-1}(\cU_x)$ is $G$-invariant and fits into the $2$-Cartesian square
    \[
    \begin{tikzcd}
    {\cW_x\simeq[V_x/H_x]} & \pt \\
    {\cU_x\simeq[V_x/\widehat G_x]} & \BG.
    \arrow[from=1-1,to=1-2]
    \arrow[from=1-1,to=2-1]
    \arrow[from=1-2,to=2-2]
    \arrow["{r|_{\cU_x}}",from=2-1,to=2-2]
    \arrow["\lrcorner"{anchor=center,pos=0.125},draw=none,from=1-1,to=2-2]
    \end{tikzcd}
    \]
    This is precisely the extension-induced quotient situation of \Cref{thm:homotopy-fixed-quotient}, with $K=H_x$ and $\widetilde G=\widehat G_x$. The normal representation can be computed in a common atlas for $\cM$ and $\cN$: since $x$ is fixed up to isomorphism, the quotient groupoid orbit through $x$ equals its original groupoid orbit, thus their normal vector spaces agree at $x$.
    
    \emph{3. The local sectors.}
    \Cref{lem:splitting-space} and \Cref{thm:homotopy-fixed-quotient} now give the sector decomposition \eqref{eq:hfp-local-model}. In particular, openness of the splitting orbits and smoothness in families are conclusions of those results. Each $V_x^\sigma$ is the intersection of $V_x$ with the fixed linear subspace of $\sigma(G)$; the centralizer is a closed compact subgroup of $H_x$. Hence the right-hand side is a differentiable stack.

    \emph{4. Evaluation on an atlas.}
    For the equivariant atlas $p:X_0\to\cM$, put $X_0^{hG}=X_0\times_{\cM}\cM^{hG}$, initially a sheaf by \Cref{rem:comparison-fibers}. For a manifold \(T\),
    \[
    X_0^{hG}(T)=\Set{(u,\sigma) |\begin{gathered}u:T\to X_0,\ \sigma:T\times G\to Y_1\\ \rs(\sigma_t(g))=\rt(\sigma_t(g))=u(t),\ \rr_1(\sigma_t(g))=g,\\
        \sigma_t(e)=1_{u(t)},\ \sigma_t(gh)=\sigma_t(g)\sigma_t(h)\end{gathered}}
    \]
    where \(\sigma_t(g)=\sigma(t,g)\). Indeed, an identification $\ell:y\to u^*p$ transports $\beta$ to $\ell\,\beta_g(g\cdot\ell)^{-1}$. The equivariance of $p$ identifies this with an arrow in the presenting groupoid, smoothly on $T\times G$ by the atlas pullback definition of $X_1$. The product law is the transported cocycle. Conversely these smooth arrows give the fixing structure on $u^*p$. Intertwiners over its identity are identities. Write $\nu:X_0^{hG}\to X_0$ for the projection.

    For an open slice $[V/L]$ of $\cN$ about a point in the image, write $\cW=[V/K]$ for its open pullback to $\cM$. Since $\cW$ is full and $G$-invariant, thus $\epsilon_{\cM}^{-1}(\cW)=\cW^{hG}$: fixing arrows of a family with values in $\cW$ belong to the restricted stack. Over $O=p^{-1}(\cW)$ put $D=O\times_{\cW}V$: the map $D\to O$ is a principal right $K$-bundle, and $D\to V$ is a submersion. With the manifold $Z=V\times_{\cW}\cW^{hG}$ of \eqref{eq:epsilon-atlas}, the slice of $W$ is
    \[
        \nu^{-1}(O)\simeq (D\times_V Z)/K,\qquad (d,z)c=(dc,c^{-1}z).
    \]
    Indeed, a lift to $D$ trivializes the underlying torsor, after which a fixing structure is a point of $Z$ over the same $v$. Descent for a different lift gives the displayed action. The fiber product is a manifold by submersiveness, and the diagonal action is free and proper since $K$ is compact and acts freely on $D$.

    These open pieces cover $X_0^{hG}$: the union of their bases $O$ is an open subset of the manifold $X_0$, so countably many suffice, and on overlaps, they represent the same sheaf, giving smooth identifications with the cocycle condition. The resulting space is manifold: points over distinct base points separate in $X_0$, and points over one base point lie in a common Hausdorff piece. This proves quasi-representability. For any manifold $T\to\cM^{hG}$,
    \[
           T\times_{\cM^{hG}}X_0^{hG}\simeq T\times_{\cM}X_0,
    \]
    thus $X_0^{hG}\to\cM^{hG}$ is a representable surjective submersion.

    \emph{5. The proper presentation.}
    The atlas groupoid is $\rX^{hG}=(X_1\times_{\rs,X_0}X_0^{hG}\rightrightarrows X_0^{hG})$. An arrow $(a,w)$ has source $w=(x,\sigma)$ and target $a w=(\rt(a),g\mapsto a\sigma(g)a^{-1})$, where $a$ has degree $e$ in $\rY$  (this is the intertwining condition above); composition is $(b,aw)(a,w)=(ba,w)$ and inversion is $(a,w)^{-1}=(a^{-1},aw)$. These maps are smooth by their action on smooth families. The source is a submersion by base change, and inversion gives the same assertion for the target. The stabilizer consists of $a\in H_x$ centralizing every $\sigma(g)$.
    
    Let $C\subset X_0^{hG}\times X_0^{hG}$ be compact with projections $C_1,C_2$. Properness of the anchor of $\rX$ makes
    \[
        B=(\rs,\rt)^{-1}\bigl(\nu(C_1)\times\nu(C_2)\bigr)\subset X_1
    \]
    compact. The preimage of $C$ under $(a,w)\mapsto(w,aw)$ is closed in $B\times C_1$: composability is a closed condition because $X_0$ is Hausdorff, and the anchor is continuous on that closed subspace. Hence the anchor of the new groupoid is proper. For a second atlas $W'$, the common pullback $X_0^{hG}\times_{\cM^{hG}}W'$ has the two principal actions by changes of frame and gives a Morita bibundle. Thus the presentation computes the intrinsic stack independently of the atlas.
\end{proof}

\subsection{Comparison over the fixed locus}
\begin{theorem}\label{thm:comparison}
    Under the hypotheses of \Cref{thm:fixed-locus-proper-regular-stack}, the morphism
    $j_{\cM}:\cM^{hG}\to\cM^G$ is a representable submersion. At a point $x$ of $\cM^G$ set $L=\widehat G_x$, $K=H_x$, $\pi=\pi_x$, $S_\pi=\operatorname{spl}(\pi)$ and $K_\sigma=C_K(\sigma)$. In the relative slice let $F=V^{L^0}$. Over the open fixed chart $[F/K]$, the comparison is equivalent to
    \[
         \coprod_{[\sigma]\in S_\pi/K}[F/K_\sigma]
              \ \longrightarrow\ [F/K].
    \]
    Each summand is a representable proper bundle with fiber $K/K_\sigma$. The local essential image of $j_{\cM}$ is the clopen substack of the fixed locus on which the extended isotropy extension splits at each point. The whole morphism $j_{\cM}$ is proper if and only if $S_{\pi_x}/H_x$ is finite at every fixed point $x$.
\end{theorem}
\begin{proof}
    A splitting has connected image containing $e$, hence $\sigma(G)\subseteq L^0$ and $F\subseteq V^{\sigma}$. Conversely, if $v$ is fixed by $\sigma(G)$, its stabilizer $L_v$ projects onto $G$. The dimension argument in \Cref{thm:fixed-locus-proper-regular-stack} then gives $v\in F$. Thus $V^{\sigma}=F$, and \eqref{eq:hfp-quotient-decomposition} identifies evaluation on the summand with the morphism induced by $K_\sigma\hookrightarrow K$ and $\id_F$. The calculation also applies when $S_\pi$ is empty.
    
    A map $T\to[F/K]$ write its object as a right $K$-torsor $P\to T$ with an equivariant map $P\to F$. A lift to $[F/K_\sigma]$ together with its frame is exactly a reduction of $P$ to $K_\sigma$. Its sheaf of choices is represented by $P/K_\sigma=P\times_K(K/K_\sigma)$. Consequently the pullback of $j_{\cM}$ is the countable disjoint union of these bundles. Each is locally trivial on $T$ with compact fiber, so it is submersive and proper. Their union is a manifold and a submersion, allowing the empty fiber. For a general $T\to\cM^G$ the open slices cover $T$, and the representing manifolds glue by uniqueness on overlaps. A countable subcover suffices and points above a common base point are separated in a single bundle chart.
    
    Every splitting over a slice center fixes all of $F$ and yields local lifts over the entire chart. Conversely a lift at $v$ yields a splitting into $L_v\subset L$. Thus the image and its complement are open; for a family this is a criterion for local lifts. On a chart with finitely many splitting classes the finite union is proper, a property local on the target for these manifold maps. Conversely, the base change at the center of a chart with infinitely many classes is an infinite disjoint union of nonempty open compact homogeneous spaces. Its open cover by summands has no finite subcover.
\end{proof}

\begin{remark}[Why a finiteness hypothesis is needed]
    For the trivial $\mathbb{S}^{1}$-action on $\mathcal B \Sone$, the extension at the unique point is $\Sone\times \Sone\to \Sone$. Its splittings are $g\mapsto(g^n,g)$, $n\in\mathbb Z$, and conjugation by $K=\Sone$ is trivial. Thus $j$ is the fold $\coprod_{n\in\mathbb Z}\mathcal B S^1\to\mathcal B \Sone$. Its pullback along the trivial torsor $\pt\to\mathcal B \Sone$ is $\mathbb Z\to\pt$ with discrete domain.
\end{remark}

\begin{remark}
    If $\cM$ is an orbifold under the same connected compact group, $K$ is finite. For two splittings the continuous map $g\mapsto\sigma_1(g)\sigma_2(g)^{-1}$ takes values in $K$ and is the identity at $e$. Connectedness makes it constant, and a splitting, when it exists, is therefore unique and every $K$-conjugate is equal to it. Thus $K_\sigma=K$, and $j$ identifies $\cM^{hG}$ with the clopen splitting substack of $\cM^G$. In particular, in this case $j$ is proper.
\end{remark}

\section{Examples}\label{sec:examples}

The distinction between $\cX^{G}$ and $\cX^{hG}$ is the distinction between existence and coherent choice.  The fixed locus is a full substack of $\cX$; the homotopy fixed stack lies over $\cX$ and may have several, one, or no points over a fixed object.

\begin{example}[Manifolds]
    For any smooth action on a manifold $M$ both stacks are the same sheaf. On a test manifold $T$ its sections are the smooth $u:T\to M$ such that $g\cdot u(t)=u(t)$ for every $g\in G$ and $t\in T$. This pointwise equality is also equality of the maps $G\times T\to M$. In the represented stack there is a unique possible fixing arrow. When the ordinary fixed set $M^G$ is embedded, this sheaf is represented by $M^G$. In particular, this holds for compact $G$.
\end{example}

\begin{example}[Two actions on the same classifying stack]\label{ex:trivial-fixed-action}
Take $G=K=\mathbb Z/2$. For the trivial action on $\mathcal B K$ the extension is $K\times G\to G$. The zero and identity homomorphisms $G\to K$ give two splittings, both centralized by $K$. Hence
\[
     (\mathcal B K)^G\simeq\mathcal B K,\qquad
     (\mathcal B K)^{hG}\simeq\mathcal B K\sqcup\mathcal B K .
\]
On the same stack consider instead the weak action induced by
$0\to K\xrightarrow{k\mapsto2k}\mathbb Z/4\to G\to0$.
Choosing representatives $s(0)=0$, $s(1)=1$ makes each individual action functor equivalent to the identity, since the extension is central. The associator is represented by the normalized cocycle
\[
 c(g,h)\equiv_2\big(s(g)+s(h)-s(g+h)\big)/2,\qquad c(1,1)=1.
\]
It represents the nonzero class in $H^2(\mathbb Z/2,\mathbb Z/2)$: a normalized coboundary has value $0$ at $(1,1)$. Local triviality of $K$-torsors and surjectivity of $\mathbb Z/4\to G$ make every object fixed. A splitting would lift the nonidentity element to an element of order two, whereas both its lifts have order four. Thus $(\mathcal B K)^G\simeq\mathcal B K$ and $(\mathcal B K)^{hG}\simeq\varnothing$ for this action.
\end{example}

\begin{example}[A proper non-regular example]
\label{ex:proper-nonregular-homotopy-fixed}
In this example the pointwise fixed set on the atlas has two branches, and the isotropy dimension jumps at their intersection.
    
    Let $K=\Sone$, $M=\mathbb C^2$, and let $K$ act diagonally on $M$. Consider the action groupoid $\rX=K\ltimes M$ and the associated differentiable stack $[M/K]$.
    The groupoid $\rX$ is proper because $K$ is compact. Its isotropy is trivial away from the origin and equals $K$ at the origin, so $\rX$ is nonregular.
    
    Let $G=\Sone$ act on $M$ by
    \[
    g\cdot(z_1,z_2)=(gz_1,g^{-1}z_2).
    \]
    The commuting actions give a strict $G$-action on $K\ltimes M$ and an induced smooth action on $[M/K]$. Equivalently, the corresponding quotient is presented by $\rY=(K\times G)\ltimes M$,
    where
    \[
        (k,g)\cdot(z_1,z_2) = (kgz_1,kg^{-1}z_2),
    \]
    and the projection $\rr:\rY\to\mathsf{B}G$ is induced on arrows by $\rr_1(k,g)=g$.
    Its fiber over the identity of $G$ is precisely $\rX=K\ltimes M$.
    
    For $z=(z_1,z_2)\in M$, the extended stabilizer is
    \begin{equation}
    \label{eq:extended-stabilizer-example}
        \widehat G_z = \Set{(k,g)\in K\times G| kgz_1=z_1,\;
           kg^{-1}z_2=z_2},
    \end{equation}
    and the canonical homomorphism to the acting group is
    \[
        \pi_z:\widehat G_z\longrightarrow G,
        \qquad
        (k,g)\longmapsto g.
    \]
    
    \emph{The fixed locus.}
    The object represented by $z$ is fixed if and only if $\pi_z$ is surjective. Suppose first that both coordinates are nonzero.
    Then the two equations in \eqref{eq:extended-stabilizer-example} give
    \[kg=1,\qquad kg^{-1}=1.\] 
    Hence $g^2=1$, so the image of $\pi_z$ is contained in the finite subgroup $\{\pm1\}\subset\Sone$. In particular, $\pi_z$ is not surjective.
    
    If $z_1\neq0$ and $z_2=0$, the choice $k=g^{-1}$ gives a lift $(g^{-1},g)$ for every $g\in G$. If $z_1=0$ and $z_2\neq0$, take $k=g$. At the origin both equations hold for every pair, so
    \[
        \widehat G_0=K\times G.
    \]
    
Consequently the set of points of the atlas pullback of $[M/K]^G$ is
\[
    (\mathbb C\times\{0\})\ \cup\ (\{0\}\times\mathbb C).
\]
The two axes meet at the origin, so their union is not an embedded submanifold of $M$. This computation describes the points of the atlas pullback. Its sections on a general test manifold must also satisfy the local lifting condition.
    
    \emph{The homotopy fixed locus.}
    The $G$- and $K$-actions commute and $K$ is abelian, so \Cref{rem:direct} gives
    \[
        [M/K]^{hG}\simeq
        \coprod_{\varphi\in\hom(G,K)}[M^\varphi/K].
    \]
    Every smooth group homomorphism $\varphi:\Sone\to\Sone$
    is of the form
    \[
        \varphi_n(g)=g^n, \qquad n\in\mathbb Z,
    \]
    so $\hom(G,K)\cong\mathbb Z$. For each integer $n$,
    \[
        M^{\varphi_n} = \Set{z\in M | (g^{n+1}z_1,g^{n-1}z_2)=(z_1,z_2)\ \forall g\in G}.
    \]
    For $n=-1$ the first character is trivial, giving $M^{\varphi_{-1}}=\mathbb C\times\{0\}$. For $n=1$ the second is trivial, giving $M^{\varphi_1}=\{0\}\times\mathbb C$. For $n\neq\pm1$ both characters are nontrivial and $M^{\varphi_n}=\{0\}$. Thus
    \[
    [\C^2/\Sone]^{h\Sone}\simeq[(\mathbb C\times\{0\})/\Sone] \sqcup [(\{0\}\times\mathbb C)/\Sone]\sqcup \coprod_{n\in\mathbb Z\setminus\{-1,1\}}\mathcal{B}\Sone.
    \]
    At the origin the distinct homotopy fixed structures are $(0,\widetilde\varphi_n)$ for $n\in\mathbb Z$. Indeed, the isotropy extension there is
    \[\begin{tikzcd}
    	1 & \Sone & {\Sone\times\Sone} & \Sone & 1
    	\arrow[from=1-1, to=1-2]
    	\arrow[from=1-2, to=1-3]
    	\arrow[from=1-3, to=1-4]
    	\arrow[from=1-4, to=1-5]
    \end{tikzcd}\]
    and every homomorphism $\varphi_n:\Sone\to \Sone$ gives a splitting
    \[
        \widetilde{\varphi}_n(g)=(g^n,g).
    \]
    The two branches of the pointwise fixed set lift to the sectors indexed by
    \[
        \varphi_{-1}(g)=g^{-1},
        \qquad
        \varphi_1(g)=g,
    \]
    respectively. They lie in distinct components of $[M/K]^{hG}$.
\end{example}

\begin{example}[The complex weighted projective line]\label{ex:weighted-proj-line}
    Let $a,b$ be positive integers. Let $K=\mathbb S^1$ act on $M=\mathbb S^3\subset\mathbb C^2$ with weights $(a,b)$ by
    \[
    k\cdot(z_1,z_2):=(k^a z_1,k^b z_2).
    \]
    The quotient stack $\mathbb P(a,b)= [M/K]$ is the \emph{differentiable complex weighted projective line}. The locally free action of the compact group $K$ makes $\mathbb P(a,b)$ an orbifold stack, so \Cref{cor:fixed-locus-orbifold-stack} and \Cref{thm:hfp-differentiability} apply. Let $G=\mathbb S^1$ act on $M$ by
    \[
    g\cdot(z_1,z_2)=(g z_1,z_2).
    \]
    This action commutes with that of $K$. It gives a strict $G$-action on $K\ltimes M$ and hence a smooth action on $\mathbb P(a,b)$.

    \emph{The fixed locus.} A point $(z_1,z_2)\in M$ represents a fixed object of $[M/K]$ if and only if, for every $g\in G$, there exists $k\in K$ such that
    \[
    (gk^a z_1,k^b z_2)=(z_1,z_2).
    \]
    If both coordinates are nonzero, $k^b=1$ and $g=k^{-a}$. Only finitely many $g$ arise in this way. If $z_2=0$, the extended stabilizer is
    \[
     \Sone\xrightarrow{\ \sim\ }\widehat G_a=\Set{(g,k)\in G\times K|g k^a=1},\qquad h\longmapsto(h^a,h^{-1}).
    \]
    Under this identification $\pi_a(h)=h^a$. If $z_1=0$, we have
    \[
    \widehat G_b=G\times\mu_b, \qquad \pi_b=\operatorname{pr}_G,
    \]
    where $\mu_c=\Set{\lambda\in\Sone|\lambda^c=1}$ is the group of $c$th roots of unity. The isotropy groups on the two axes are $\mu_a$ and $\mu_b$, respectively. The inverse image in $M$ of the pointwise fixed locus is therefore
    \[
    F=\Set{(z_1,0)| \,|z_1|=1}\cup\Set{(0,z_2)|\, |z_2|=1},
    \]
    and the fixed locus stack is
    \[
    \mathbb P(a,b)^G\simeq \mathcal{B}\mu_a\sqcup \mathcal{B}\mu_b.
    \]
    \emph{The homotopy fixed locus.} As in the preceding example, the characters $\varphi_n(g)=g^n$ index the sectors in \Cref{rem:direct}. For $n\in\mathbb Z$,
    \[
        M^{\varphi_n} = \Set{z\in M | (g^{an+1}z_1,g^{bn}z_2)=(z_1,z_2)\ \forall g\in G}.
    \]
    When $a=1$, only $n=-1$ and $n=0$ contribute. For $n=-1$ we have
    \[
    M^{\varphi_{-1}}=\Set{(z_1,0) |\, |z_1|=1},\qquad [M^{\varphi_{-1}}/K]\simeq \pt.
    \]
    If $n=0$, then
    \[
    M^{\varphi_0}=\Set{(0,z_2)|\, |z_2|=1}, \qquad [M^{\varphi_0}/K]\simeq\mathcal{B}\mu_b.
    \]
    For $a>1$, the equation $an+1=0$ has no integer solution, so only $n=0$ contributes. Combining the two cases gives
    \[
     \mathbb P(a,b)^{hG}\simeq
     \begin{cases}
     \mathcal{B}\mu_a\sqcup \mathcal{B}\mu_b, &a=1,\\
     \mathcal{B}\mu_b, &a>1.
     \end{cases}
    \]
\end{example}

\begin{example}[A disconnected acting group]\label{ex:disconnected-cross}
    Let $K=\{1,-1\}$ act on $\mathbb R^2$ by simultaneous sign change, and let $G=\mathbb Z/2$ act by $(a,b)\mapsto(a,-b)$. These commuting actions give a strict $G$-action on the proper \'etale groupoid $K\ltimes\mathbb R^2$ presenting $\cM=[\mathbb R^2/K]$. The extended stabilizer maps onto $G$ exactly when some $k\in K$ satisfies $k(a,-b)=(a,b)$. The choice $k=1$ gives $b=0$ and $k=-1$ gives $a=0$. This is the cross in \cite[Example~4.10]{jiang2026equivariant}. The isotropy groups of $\cM$ have dimension zero everywhere.
    
    In families over the atlas, lifting the constant section with value the generator of $G$ amounts to choosing a locally constant sign. Conversely, such a choice gives lifts for any smooth $G$-valued section since $G$ is discrete. Thus a fixed family is locally contained in one axis, and the pointwise cross cannot be the image of an embedded fixed submanifold. The two homomorphisms $G\to K$ give
    \[
       \cM^{hG}\simeq
       [(\mathbb R\times\{0\})/K]\ \sqcup\ [(\{0\}\times\mathbb R)/K].
    \]
    At the origin there are two distinct fixing structures, each with isotropy $K$. Thus the embeddedness assertion of \Cref{thm:fixed-locus-proper-regular-stack} can fail for a disconnected group, even when the isotropy dimension is constant. \Cref{thm:hfp-differentiability} still applies to the homotopy fixed locus.
\end{example}

\appendix
\section{Kan fibrations and stack actions}\label[appendix]{appendix:translation}

Jiang, Yang and Chen introduce the Kan presentation, its bibundles $P_g$ and the free $G$-groupoid $\mathsf Z$ with $\mathsf Z/G\cong\rY$ in \cite[\S~2.2]{jiang2026equivariant}. We identify the corresponding smooth action and its quotient as stacks. The division maps below also control smooth dependence on the group parameters.

\begin{definition}[{\cite[Definition 2.4]{jiang2026equivariant}}]\label{def:kan}
    A homomorphism of Lie groupoids $\rr:\rY\to\mathsf{B}G$ is a \emph{Kan fibration} if
    \[
    (\rs_{\rY},\rr_1):Y_1\longrightarrow Y_0\times G
    \]
    is a surjective submersion. Equivalently, $(\rt_{\rY},\rr_1)$ is a surjective submersion. Its kernel is the Lie groupoid
    \[
    \rX=\ker(\rr)=\bigl(\rr_1^{-1}(e)\rightrightarrows X_0\bigr),
    \qquad X_0=Y_0.
    \]
\end{definition}

The equivalence of the two conditions follows by inversion in $\rY$ and $G$. In particular $\rr_1$ is a submersion, and the source and target maps of $\rX$ are surjective submersions by base change along $Y_0\times\{e\}$. Set
\[
\cM=[X_0/X_1],\qquad \cN=[Y_0/Y_1],\qquad
r:\cN\longrightarrow\BG,
\]
Thus $r$ is the morphism of classifying stacks obtained by applying the stack construction to the groupoid functor $\rr$.

For $g\in G$, put $P_g=\rr_1^{-1}(g)$. Its left and right moment maps are, respectively, $\rt_{\rY}$ and $\rs_{\rY}$, and the commuting $\rX$-actions are
\[
c\cdot a=ca,\qquad a\cdot d=ad,
\qquad c,d\in X_1,\quad a\in P_g,
\]
whenever the arrows are composable. We use the bibundle convention in which a morphism from a source groupoid to a target groupoid is represented by a left-principal target--source bibundle: the left target action is principal over the right anchor map to the source object manifold.

\begin{lemma}\label{lem:Pg-principal}
    For every $g\in G$, $P_g$ is a Morita self-equivalence bibundle of $\rX$. Multiplication in $\rY$ induces bibundle isomorphisms
    \[
    c_{g,h}:P_g\otimes_{\rX}P_h\xrightarrow{\sim}P_{gh},
    \qquad[a,b]\longmapsto ab,
    \]
    and $P_e=X_1$ is the identity bibundle. These isomorphisms depend smoothly on $(g,h)$.
\end{lemma}
\begin{proof}
    Both anchors on $P_g$ are surjective submersions by base change of the Kan map and its target version. For $a,b\in P_g$ with the same source the unique $c$ such that $ca=b$ is $c=ba^{-1}\in X_1$. For $a,b$ with the same target the unique $d$ such that $ad=b$ is $d=a^{-1}b\in X_1$. These smooth division maps are inverses to the two principality maps.

    Multiplication identifies $(ad,b)$ with $(a,db)$ in the balanced product. Given $z\in P_{gh}$, choose $b\in P_h$ with $\rs(b)=\rs(z)$ and put $a=zb^{-1}\in P_g$. A second choice $b'$ gives $d=b'b^{-1}\in X_1$ and replaces $(a,b)$ by $(ad^{-1},db)$. Thus $z\mapsto[zb^{-1},b]$ is inverse to multiplication. For varying $(g,h,z)$, pull local smooth sections of $(\rs_{\rY},\rr_1)$ back along $(g,h,z)\mapsto(\rs(z),h)$. They give smooth local formulas for this inverse and agree in the balanced quotient. Multiplication is smooth before passing to the quotient, so the isomorphisms are smooth in the parameters.
\end{proof}

The individual equivalences represented by $P_g$ assemble into a smooth action morphism. Indeed, the whole manifold $Y_1$ is a left-principal bibundle from $\underline G\times\rX$ to $\rX$, where $\underline G=(G\rightrightarrows G)$ is the unit groupoid. Its left anchor map is $\rt_{\rY}$ and its right anchor map is
\[
(\rr_1,\rs_{\rY}):Y_1\longrightarrow G\times X_0.
\]
The left action is postcomposition by $X_1$ and the right action is precomposition by $X_1$ at the fixed parameter $g$. The same division formula as in the lemma proves left principality. Hence this bibundle presents
\[
\mu:G\times\cM\longrightarrow\cM.
\]
Here $\underline G\times\rX$ presents $G\times\cM$. Multiplication supplies the associator $\alpha_{g,h}:\mu_g\mu_h\Rightarrow\mu_{gh}$, and the identity bibundle supplies the unitor $\iota:\id_{\cM}\Rightarrow\mu_e$.

\begin{theorem}[From a Kan fibration to a stack action]\label{thm:kan-to-action}
    A Kan fibration $\rr:\rY\to\mathsf{B}G$ determines, up to coherent equivalence, a smooth $G$-action on $\cM=[X_0/X_1]$. The action of $g$ is represented by $P_g$, with associator and unitor as above. There is an equivalence
    \[
    [\cM/G]\simeq\cN
    \]
    compatible with the morphisms to $\BG$. Under this equivalence, $r$ classifies the principal $G$-bundle $q:\cM\to\cN$ induced by $\rX\hookrightarrow\rY$, and the square
    \begin{equation}\label{eq:appendix-torsor-square}
    \begin{tikzcd}
    \cM & \pt \\
    \cN & \BG
    \arrow[from=1-1,to=1-2]
    \arrow["q"',from=1-1,to=2-1]
    \arrow[from=1-2,to=2-2]
    \arrow["r",from=2-1,to=2-2]
    \arrow["\lrcorner"{anchor=center,pos=0.125},draw=none,from=1-1,to=2-2]
    \end{tikzcd}
    \end{equation}
    is $2$-Cartesian, with the $2$-cell supplied by the trivialization on the kernel.
\end{theorem}
\begin{proof}
    Multiplication in the parameterized bibundles sends both bracketings of $[a,b,c]$ to $abc$, and composition with a unit gives the triangle equations. To identify the principal bundle, use the ordinary pullback of groupoids $\mathsf Z=\rY\times_{\mathsf B G}\mathsf E G$ of \cite[\S 2.2]{jiang2026equivariant}, writing $ba$ for $b$ after $a$. Here $\mathsf E G=(G\times G\rightrightarrows G)$ is the pair groupoid of $G$ and $(b,a):a\to b$ has degree $ba^{-1}$. Then
    \[
     Z_0=Y_0\times G,\qquad
     (a,k):(\rs(a),k)\longrightarrow(\rt(a),\rr_1(a)k),\qquad
     (b,\rr_1(a)k)(a,k)=(ba,k).
    \]
    Right translation in the second coordinate is free and gives $\mathsf Z/G\cong\rY$ on objects and arrows. The inclusion of the full restriction at $k=e$ identifies it with $\rX$. Its essential-surjectivity map $a\mapsto(\rt(a),\rr_1(a))$ is a surjective submersion. Hence this inclusion is a Morita equivalence.
    
    On the atlas $Y_0$ the torsor classified by $r$ is $Y_0\times G$ with descent $a:(\rs(a),k)\mapsto(\rt(a),\rr_1(a)k)$. Its total stack is presented by $\mathsf Z$. This identifies $[Z_0/Z_1]$ with $\cN\times_{\BG}\pt$ on all test manifolds: a lift is a trivialization of the pulled-back torsor, and changes of atlas lift act by the displayed descent formula. The frame $k=e$ on kernel arrows gives the $2$-cell in \eqref{eq:appendix-torsor-square}. Written as a left action, translation sends $(x,e)$ to $(x,g^{-1})$. An arrow $(x,g^{-1})\to(y,e)$ is exactly an element of $P_g$ from $x$ to $y$. Its multiplication is $c_{g,h}$, identifying this action equivariantly with the bibundle action. Taking the quotient gives $[\cM/G]\simeq\cN$ over $\BG$.
\end{proof}

Conversely, let $(\cM,\mu,\iota,\alpha)$ be a differentiable $G$-stack. Choose a $G$-atlas $a:X_0\to\cM$ as in \Cref{prop:equivariant-atlas}, and let $\rX=(X_1\rightrightarrows X_0)$ be its induced strict equivariant presentation. The composite $X_0\to[\cM/G]$ is an atlas; compare \cite[Proposition~2.12]{barbosa2021equivariant}. Its groupoid is the one used in the proof of \Cref{prop:proper-quotient}:
\[
Y_0=X_0,\qquad Y_1=G\times X_1,
\]
\[
\rs_{\rY}(g,\gamma)=g^{-1}\cdot \rs_{\rX}(\gamma),\qquad
\rt_{\rY}(g,\gamma)=\rt_{\rX}(\gamma).
\]
An arrow $(g,\gamma):x\to y$ means $\gamma:g\cdot x\to y$ in $\rX$, composition is
\[
(h,\delta)\circ(g,\gamma)
=(hg,\delta\circ(h\cdot\gamma)),
\]
and define $\rr_1(g,\gamma)=g$.

\begin{theorem}[From a stack action to a Kan fibration]\label{thm:action-to-kan}
    The homomorphism $\rr:\rY\to\mathsf{B}G$ just constructed is a Kan fibration. Its kernel is canonically isomorphic to $\rX$, and the action reconstructed by \Cref{thm:kan-to-action} is equivariantly equivalent to the original action on $\cM$.
\end{theorem}
\begin{proof}
    The composition formula shows that $\rr$ is a groupoid homomorphism. The map $(\rs_{\rY},\rr_1)$ is the composite
    \[
    G\times X_1\xrightarrow{\id_G\times \rs_{\rX}}G\times X_0
    \xrightarrow{\ d\ }X_0\times G,
    \qquad d(g,z)=(g^{-1}\cdot z,g).
    \]
    The first map is a surjective submersion and $d$ is a diffeomorphism, with inverse $(x,g)\mapsto(g,g\cdot x)$. Hence $\rr$ satisfies the Kan condition at every arrow.

    Its kernel has arrows $\{e\}\times X_1$, with the original structure maps of $\rX$. The fiber $\rr_1^{-1}(g)$ consists exactly of the arrows $\gamma:g\cdot x\to y$, so its bibundle represents the original action of $g$ through the equivariance $2$-cell of $a$. The composition of these arrows recovers the associator, and the identity arrows recover the unitor $\iota$. Since the atlas equivariance satisfies the action coherence identities, these identifications give an equivariant equivalence of the reconstructed action with $(\mu,\iota,\alpha)$.
\end{proof}

\begin{corollary}\label{cor:equivalence}
    Smooth $G$-actions on differentiable stacks, up to equivariant equivalence, correspond to Kan presentations up to Morita equivalence over $\BG$, with a specified identification of the kernel stack. This correspondence is independent of the chosen equivariant atlas.
\end{corollary}
\begin{proof}
It follows directly from \Cref{thm:kan-to-action,thm:action-to-kan}.
\end{proof}

\begin{remark}
    This comparison also identifies the stabilizers used in the main text. For $x\in X_0=Y_0$, the extended stabilizer is the isotropy group $\Aut_{\rY}(x)$, with projection $\rr_1|_{\Aut_{\rY}(x)}$ to $G$ and kernel $\Aut_{\rX}(x)$. Surjectivity of this projection is the pointwise fixed-locus condition, while a homotopy fixed structure is a chosen smooth homomorphic section. The family condition remains the sheaf-theoretic condition of \Cref{def:fixed-locus}.
\end{remark}

\begin{remark}[The categorical interpretation of the correspondence]\label{rem:quotient-slice-categorical}
The preceding constructions give a biequivalence
\[
   \mathbf{DiffSt}_G\simeq\mathbf{DiffSt}/\BG,\qquad
   \cM\longmapsto([\cM/G],r),\qquad
   (\cN,r)\longmapsto\cN\times_{\BG}\pt .
\]
Indeed $r$ classifies a principal $G$-bundle on $\cN$ with its own pullback on every test manifold \cite[Propositions~3.2--3.3]{biswas2023atiyah}. Its total stack is differentiable: over an atlas it is an ordinary principal bundle, and its total space provides an atlas by descent. In the bicategorical slice a morphism $(f,\theta)$ induces the equivariant map on frames just described. A $2$-cell $\eta:f\Rightarrow f'$ is allowed precisely when $\theta'\circ(r'\eta)=\theta$; this is precisely compatibility with those frames. In the other direction an equivariant map descends to a map on the quotients with the induced identification of their torsors, and an equivariant $2$-cell descends along the quotient atlas. Local trivializations of a torsor give the unit and counit equivalences; a change of trivialization acts by its transition function, so the equivalences commute with pullback and composition.

In this description $\Trivial(\cM)$ corresponds to $\BG\times\cM\to\BG$, $\Forgetful$ is base change along $\pt\to\BG$, and \eqref{eq:hfp-section-stack} computes the stack of sections. If $M$ is a manifold, passing to isomorphism classes gives $\pi_0(\Mor_{\St}(M,\BG))\cong[M,BG]$, where $BG$ is the topological classifying space and the right side denotes ordinary homotopy classes.
\end{remark}

\clearpage

\printbibliography

@misc{alper2025stacks,
 author={Alper, Jarod}, 
 title={Stacks and Moduli},
 year={2026}, 
 note={Manuscript, January 5, 2026},
 url={https://sites.math.washington.edu/~jarod/moduli-versions/moduli-1-5-26.pdf}}

@article{aranha2025virtual,
 author={Aranha, Dhyan and Khan, Adeel A. and Latyntsev, Alexei and Park, Hyeonjun and Ravi, Charanya},
 title={Virtual localization revisited}, 
 journal={Advances in Mathematics},
 volume={479}, 
 year={2025}, 
 pages={110434},
 doi={10.1016/j.aim.2025.110434}}

@article{barbosa2021equivariant,
 author={Barbosa-Torres, Luis Alejandro and Neumann, Frank}, 
 title={Equivariant cohomology for differentiable stacks},
 journal={Journal of Geometry and Physics}, 
 volume={160}, 
 year={2021}, 
 pages={103974}, 
 doi={10.1016/j.geomphys.2020.103974}}

@article{behrend2011differentiable,
 author={Behrend, Kai and Xu, Ping}, 
 title={Differentiable stacks and gerbes},
 journal={Journal of Symplectic Geometry}, 
 volume={9}, 
 year={2011}, 
 pages={285--341},
 doi={10.4310/JSG.2011.v9.n3.a2}}

@article{biswas2023atiyah,
 author={Biswas, Indranil and Chatterjee, Saikat and Koushik, Praphulla and Neumann, Frank},
 title={Atiyah sequences and connections on principal bundles over {Lie} groupoids and differentiable stacks},
 journal={Journal of Noncommutative Geometry}, 
 volume={17}, 
 year={2023}, 
 pages={407--437}, 
 doi={10.4171/JNCG/486}}

@article{bursztyn2020principal,
 author={Bursztyn, Henrique and Noseda, Francesco and Zhu, Chenchang}, 
 title={Principal actions of stacky {Lie} groupoids},
 journal={International Mathematics Research Notices}, 
 year={2020}, 
 volume={2020}, 
 number={16}, 
 pages={5055--5125}, 
 doi={10.1093/imrn/rny142}}

@article{cardenas2020stability,
 author={C{\'a}rdenas, Cristian Camilo and Struchiner, Ivan}, 
 title={Stability of {Lie} group homomorphisms and {Lie} subgroups},
 journal={Journal of Pure and Applied Algebra}, 
 volume={224}, 
 year={2020}, 
 pages={1280--1296},
 doi={10.1016/j.jpaa.2019.07.017}}

@article{carlsson1992survey,
 author={Carlsson, Gunnar}, 
 title={A survey of equivariant stable homotopy theory}, 
 journal={Topology}, 
 volume={31}, 
 year={1992}, 
 pages={1--27},
 doi={10.1016/0040-9383(92)90061-L}}

@article{chenruan,
 author={Chen, Weimin and Ruan, Yongbin}, 
 title={A new cohomology theory of orbifold},
 journal={Communications in Mathematical Physics}, 
 volume={248}, 
 year={2004}, 
 pages={1--31}, 
 doi={10.1007/s00220-004-1089-4}}

@phdthesis{cnossen2024twisted,
 author={Cnossen, Bastiaan}, 
 title={Twisted ambidexterity in equivariant homotopy theory: two approaches},
 school={Universit{\"a}t Bonn}, 
 year={2024},
 url={https://bastiaancnossen.com/thesis.pdf}}

@book{conner1962differentiable,
 author={Conner, P. E. and Floyd, E. E.}, 
 title={Differentiable Periodic Maps},
 publisher={Springer}, 
 year={1964}, 
 series={Ergebnisse der Mathematik und ihrer Grenzgebiete}, 
 volume={33},
 doi={10.1007/978-3-662-41633-4}}

@article{crainic2013linearization,
 author={Crainic, Marius and Struchiner, Ivan}, 
 title={On the linearization theorem for proper {Lie} groupoids},
 journal={Annales scientifiques de l'{\'E}cole Normale Sup{\'e}rieure}, 
 volume={46}, 
 year={2013}, 
 pages={723--746}, 
 doi={10.24033/asens.2200}}

@misc{del2012lie,
 author={del Hoyo, Mat{\'i}as L.}, 
 title={{Lie} groupoids and differentiable stacks},
 year={2012}, 
 note={arXiv:1212.6714. Published as \emph{Lie groupoids and their orbispaces} in \emph{Portugaliae Mathematica} 70 (2013), 161--209},
 doi={10.4171/PM/1930}}

@misc{ginot2013introduction,
 author={Ginot, Gr{\'e}gory}, title={Introduction to Differentiable Stacks (and Gerbes, Moduli Spaces \ldots)},
 year={2013}, 
 note={Lecture notes, Galileo Galilei Institute, Florence},
 url={https://www.math.univ-paris13.fr/~ginot/papers/DiffStacksIGG2013.pdf}}

@misc{ginotnoohi,
 author={Ginot, Gr{\'e}gory and Noohi, Behrang}, title={Group actions on stacks and applications to equivariant string topology for stacks},
 year={2015},
 %url={https://arxiv.org/abs/1206.5603},
 doi={10.48550/arXiv.1206.5603}}

@incollection{heinloth2005notes,
 author={Heinloth, Jochen}, title={Notes on differentiable stacks},
 booktitle={Mathematisches Institut, Seminars, 2004--05}, publisher={Universit{\"a}t G{\"o}ttingen}, 
 year={2005}, 
 pages={1--32},
 %note={DOI of the edited volume \url{https://doi.org/10.17875/gup2005-48}},
 doi={10.17875/gup2005-48}}

@misc{jiang2026equivariant,
 author={Jiang, Fengyu and Yang, Yang and Chen, Bohui}, 
 title={Equivariant basic cohomology of {Lie} groupoids},
 year={2026},
% url={https://arxiv.org/abs/2606.06163v1},
 doi={10.48550/arXiv.2606.06163}}

@article{lee1969conjugacy,
 author={Lee, Dong Hoon and Wu, Ta-Sun}, 
 title={On conjugacy of homomorphisms of topological groups},
 journal={Illinois Journal of Mathematics}, 
 volume={13}, 
 year={1969}, 
 pages={694--699},
 doi={10.1215/ijm/1256053429}}

@article{lee1974deformations,
 author={Lee, Dong Hoon}, 
 title={On deformations of homomorphisms of locally compact groups},
 journal={Transactions of the American Mathematical Society}, volume={191}, 
 year={1974}, 
 pages={353--361}, 
 doi={10.1090/S0002-9947-1974-0435284-6}}

@article{lerman2012hamiltonian,
 author={Lerman, Eugene and Malkin, Anton}, title={Hamiltonian group actions on symplectic {Deligne--Mumford} stacks and toric orbifolds},
 journal={Advances in Mathematics}, 
 volume={229}, 
 year={2012}, 
 pages={984--1000}, 
 doi={10.1016/j.aim.2011.10.013}}

@misc{moerdijk2002classification,
 author={Moerdijk, Ieke}, title={On the classification of regular {Lie} groupoids}, year={2002}, note={arXiv:math/0203099v1},
 doi={10.48550/arXiv.math/0203099}}

@book{montgomery1940topological,
 author={Montgomery, Deane and Zippin, Leo}, 
 title={Topological Transformation Groups},
 publisher={Interscience}, 
 year={1955}}

@article{Mrcun2018,
 author={Mr{\v c}un, Janez}, 
 title={Monodromy and faithful representability of {Lie} groupoids},
 journal={Topology and its Applications}, 
 volume={240}, 
 year={2018}, 
 pages={137--158}, 
 doi={10.1016/j.topol.2018.03.011}}

@article{nijenhuis1967deformations,
 author={Nijenhuis, Albert and Richardson, R. W.}, 
 title={Deformations of homomorphisms of {Lie} groups and {Lie} algebras},
 journal={Bulletin of the American Mathematical Society}, 
 year={1967}, 
 volume={73}, 
 pages={175--179},
 doi={10.1090/S0002-9904-1967-11703-8}}

@article{rezk2018classifying,
 author={Rezk, Charles}, 
 title={Classifying spaces for 1-truncated compact {Lie} groups},
 journal={Algebraic \& Geometric Topology}, 
 volume={18}, 
 year={2018}, 
 pages={525--546}, 
 doi={10.2140/agt.2018.18.525}}

@book{riehl2017category,
 author={Riehl, Emily}, 
 title={Category Theory in Context}, 
 publisher={Dover Publications}, 
 year={2016},
 url={https://emilyriehl.github.io/files/context-v1.pdf}}

@article{romagny2005group,
 author={Romagny, Matthieu}, 
 title={Group actions on stacks and applications},
 journal={Michigan Mathematical Journal}, 
 volume={53}, 
 year={2005}, 
 pages={209--236},
 doi={10.1307/mmj/1114021093}}

@misc{romagny2022algebraicity,
 author={Romagny, Matthieu}, 
 title={Algebraicity and smoothness of fixed point stacks},
 year={2022},
% url={https://arxiv.org/abs/2205.11114v3},
 doi={10.48550/arXiv.2205.11114}}

@book{schwede,
 author={Schwede, Stefan}, 
 title={Global Homotopy Theory}, 
 year={2018},
 publisher={Cambridge University Press}, 
 series={New Mathematical Monographs}, 
 volume={34}, 
 doi={10.1017/9781108349161}}

@incollection{thomason1983homotopy,
 author={Thomason, R. W.}, 
 title={The homotopy limit problem}, booktitle={Proceedings of the Northwestern Homotopy Theory Conference},
 series={Contemporary Mathematics}, 
 volume={19}, 
 publisher={American Mathematical Society}, 
 year={1983}, 
 pages={407--419},
 doi={10.1090/conm/019/711065}}

@article{weil1964remarks,
 author={Weil, Andr{\'e}}, 
 title={Remarks on the cohomology of groups},
 journal={Annals of Mathematics}, 
 year={1964}, 
 volume={80}, 
 pages={149--157},
 doi={10.2307/1970495}}

@article{zhu2009n,
 author={Zhu, Chenchang}, 
 title={$n$-groupoids and stacky groupoids}, 
 journal={International Mathematics Research Notices},
 year={2009},
 volume={2009}, 
 number={21}, 
 pages={4087--4141}, 
 doi={10.1093/imrn/rnp080}}

\end{document}